\documentclass{article}

\usepackage[pdftex,
pdfauthor={Anirban Chaudhuri, Boris Kramer, Matthew Norton, Johannes O. Royset, Karen Willcox},
pdftitle={Certifiable Risk-Based Engineering Design Optimization},
pdfkeywords={CVaR, buffered probability of exceedence, failure probability, risk-averse optimization, uncertainty quantification, reliability-based design optimization},
pdfproducer={Latex with hyperref},
pdfcreator={pdflatex}]{hyperref}

\usepackage[margin=1in]{geometry} 
\usepackage{graphicx}

\usepackage{wrapfig}
\usepackage{threeparttable}
\usepackage{multirow}
\usepackage{dcolumn}
 \newcolumntype{d}{D{.}{.}{-1}}
\usepackage{nomencl}
 \makeglossary
\usepackage{subfigure}
\usepackage{subfigmat}
\usepackage{fancyvrb}
 \fvset{fontsize=\footnotesize,xleftmargin=2em}
\usepackage{enumitem}

\usepackage{amsfonts}
\usepackage{amsmath}
\usepackage{amssymb}
\usepackage{cleveref}
\usepackage{mathtools}

\newtheorem{remark}{Remark}

\usepackage{algorithm,algpseudocode}

\usepackage{color}
\usepackage{tikz}
\usetikzlibrary{shapes,arrows,shapes.multipart,calc,backgrounds, positioning, fit, decorations.pathreplacing,decorations.markings}
\usepackage{pgfplots}

\title{Certifiable Risk-Based Engineering Design Optimization}

\author{%
	Anirban Chaudhuri\footnote{Research Scientist, Department of Aeronautics and Astronautics, anirbanc@mit.edu.} \\ 
	{\normalsize\itshape Massachusetts Institute of Technology, Cambridge, MA, 02139, USA} \smallskip \\ 
	Boris Kramer\footnote{Assistant Professor, Department of Mechanical and Aerospace Engineering, bmkramer@ucsd.edu.} \\
	{\normalsize\itshape University of California San Diego, CA, 92093, USA} \smallskip \\
	Matthew Norton\footnote{Assistant Professor, Department of Operations Research, mdnorto@gmail.com.},
	\ 
	Johannes O. Royset\footnote{Professor, Department of Operations Research, joroyset@nps.edu.} \\
	{\normalsize\itshape Naval Postgraduate School, Monterey, CA, 93943, USA} \smallskip\\
	Karen E. Willcox\footnote{Director, Oden Institute for Computational Engineering and Sciences, kwillcox@oden.utexas.edu}  \\ {\normalsize\itshape University of Texas at Austin, Austin, TX, 78712, USA}
}
\date{}

\newcommand{\failset}{\mathcal{G}}
\newcommand{\bz}{\boldsymbol{z}}
\newcommand{\bd}{\boldsymbol{d}}
\newcommand{\bb}{\boldsymbol{b}}

\newcommand{\alphaT}{\alpha_{\scriptscriptstyle \text{T}}}
\newcommand{\alphaTrem}{\alpha_{\scriptscriptstyle \text{\emph{T}}}}

\newcommand{\qbar}{\overline{Q}}
\newcommand{\pt}{p_t}
\newcommand{\hpt}{\hat{p}_t}
\newcommand{\pbuf}{\overline{p}_t}

\newcolumntype{L}[1]{>{\raggedright\let\newline\\\arraybackslash\hspace{0pt}}m{#1}}
\newcolumntype{C}[1]{>{\centering\let\newline\\\arraybackslash\hspace{0pt}}m{#1}}

\graphicspath{{figures_riskOUU/}}

\begin{document}
	
\pgfmathdeclarefunction{frechet}{3}{%
	\pgfmathparse{(#1/#2)*((x-#3)/#2)^(-1-#1)*exp(-((x-#3)/#2)^(-#1))}%
}

\maketitle

\begin{abstract} 
Reliable, risk-averse design of complex engineering systems with optimized performance requires dealing with uncertainties. A conventional approach is to add safety margins to a design that was obtained from deterministic optimization. Safer engineering designs require appropriate cost and constraint function definitions that capture the \textit{risk} associated with unwanted system behavior in the presence of uncertainties. The paper proposes two notions of certifiability. The first is based on accounting for the magnitude of failure to ensure data-informed conservativeness. The second is the ability to provide optimization convergence guarantees by preserving convexity. Satisfying these notions leads to \textit{certifiable} risk-based design optimization (CRiBDO). In the context of CRiBDO, risk measures based on superquantile (a.k.a.\ conditional value-at-risk) and buffered probability of failure are analyzed. CRiBDO is contrasted with reliability-based design optimization (RBDO), where uncertainties are accounted for via the probability of failure, through a structural and a thermal design problem. A reformulation of the short column structural design problem leading to a convex CRiBDO problem is presented. The CRiBDO formulations capture more information about the problem to assign the appropriate conservativeness, exhibit superior optimization convergence by preserving properties of underlying functions, and alleviate the adverse effects of choosing hard failure thresholds required in RBDO.
\end{abstract}

\section{Introduction} \label{s:introduction}
The design of complex engineering systems requires quantifying and accounting for risk in the presence of uncertainties. This is not only vital to ensure safety of designs but also to safeguard against costly design alterations late in the design cycle. The traditional approach is to add safety margins to compensate for uncertainties \textit{after} a deterministic optimization is performed. This produces a sense of security, but is at best an imprecise recognition of risk and results in overly conservative designs that can limit performance. Properly accounting for risk \textit{during} the design optimization of those systems could allow for more efficient designs. For example, payload increases for spacecraft and aircraft could be possible \textit{without sacrificing safety}. The financial community has long recognized the superiority of specific risk measures in portfolio optimization (most importantly the conditional-value-at-risk (CVaR) pioneered by Rockafellar and Uryasev~\cite{RTRockafellar_SUryasev_2002a}), see~\cite{RTRockafellar_SUryasev_2002a, PKrokhmal_JPalmquist_SUryasev_2002a, RMansini_WOgryczak_MGSperanza_2007a}. In the financial context, it is understood that exposure to tail risk---rather rare events---can lead to catastrophic outcomes for companies, and adding too many ``safety factors'' (insurance, hedging) reduces profit. Analogously, in the engineering context, the problem is to find safe engineering designs without unnecessarily limiting performance and limiting the effects of the heuristic guesswork of choosing thresholds. 

In general, there are two main issues when formulating a design optimization under uncertainty problem: (1) what to optimize and (2) how to optimize. The first issue involves deciding the design criterion, which in the context of decision theory could boil down to what type of utility function to use. What is a meaningful way of making design decisions under uncertainty? One would like to have a framework that can reflect stakeholders' preferences, but at the same time is relatively simple and can be explained to the public, to a governor, to a CEO, etc. The answer for what to optimize directly influences how you optimize. If the ``what to optimize'' was chosen poorly, the second issue becomes much more challenging. Design optimization of a real-world system is difficult, even in a deterministic setting, so it is essential to manage complexity as we formulate the design-under-uncertainty problem. Thus, any design criterion that preserves convexity and other desirable mathematical properties of the underlying functions is preferable as it simplifies the subsequent optimization.

This motivates us to incorporate specific mathematical measures of risk, either as a design constraint or cost function, into the design optimization formulation. To this end, we focus on two particular risk measures that have potentially superior properties: (i) superquantile/CVaR~\cite{rockafellar2013superquantiles,rockafellar2000optimization}, and (ii) buffered probability of failure (bPoF)~\cite{rockafellar2010buffered}. Three immediate benefits of using these risk measures arise. First, both risk measures recognize extreme (tail) events which automatically enhances resilience. Second, they preserve convexity of underlying functions so that specialized and provably convergent optimizers can be employed. This drastically improves optimization performance.  Third, superquantile and bPoF are conservative risk measures that add a buffer zone to the limiting threshold by taking into account the magnitude of failure. This can be handled by adding safety factors to the threshold; however, it has been shown before that probabilistic approaches lead to safer designs with optimized performance compared to the safety factor approach~\cite{roland1990system,moller2008principles,suzuki2014analytical}.
Superquantile/CVaR has been recently used in specific formulations in  civil~\cite{RTRockafellar_JORoyset_2015a,WZhang_HRahimian_GBayraksan_2016a}, naval~\cite{JORoyset_LBonfiglio_GVernengo_SBrizzolara_2017a,bonfiglio2019multidisciplinary} and aerospace~\cite{HYang_MGunzburger_2016a,chaudhuri2020multifidelity} engineering, as well as general
PDE-constrained optimization~\cite{DPKouri_TMSurowiec_2016a,ZZou_DPKouri_WAquino_2018a,HKTQ18CVaRROMS,HKT2020_Adaptive_ROM_CVAR_estimation}. The bPoF risk measure has been shown to possess beneficial properties when used in optimization~\cite{rockafellar2010buffered,mafusalov2018buffered,norton2019calculating,rockafellar2020minimizing}, yet has been seldom used in engineering to-date~\cite{basova2011computational,minguez2013iterative,MMHarajli_RTRockafellar_JORoyset_2015a,royset2019risk}. We contrast these above risk-based engineering design methods with the most common approach to address parametric uncertainties in engineering design, namely reliability-based design optimization (RBDO)~\cite{yao2011review,aoues2010benchmark} which uses the probability of failure (PoF) as a design constraint. We discuss the specific advantages of using these ways of measuring risk in the design optimization cycle and their effect on the final design under uncertainty.

In this paper, we define two certifiability conditions for risk-based design optimization that can certify designs against near-failure and catastrophic failure events, and guarantee convergence to the global optimum based on preservation of convexity by the risk measures. We call the optimization formulations using risk measures satisfying any of the certifiability conditions as \textbf{C}ertifiable \textbf{Ri}sk-\textbf{B}ased \textbf{D}esign \textbf{O}ptimization (CRiBDO). Risk measures satisfying both certifiability conditions lead to strongly certifiable risk-based design. We analyze superquantile and bPoF, which are examples of risk measures satisfying the certifiability conditions. We discuss how the nature of probabilistic conservativeness introduced through superquantile and bPoF makes practical sense since it is data-informed and based on the magnitude of failure. The data-informed probabilistic conservativeness of superquantiles and bPoF circumvents the guesswork associated with setting safety factors (especially, for the conceptual design phase) and transcends the limitations of setting hard thresholds for limit state functions used in PoF. This helps us move away from \textit{being conservative blindly to being conservative to the level dictated by the data}. We compare the different risk-based design optimization formulations using a structural and a thermal design problem. For the structural design of a short column problem, we show a convex reformulation of the objective and limit state functions that leads to a convex CRiBDO formulation. 

The remainder of this paper is organized as follows. We summarize the widely-used RBDO formulation in Section~\ref{s:2}. The different risk-based optimization problem formulations along with the risk measures used in this work are described in Section~\ref{s:3}. Section~\ref{s:4} explains the features of different risk-based optimization formulations through numerical experiments on the short column problem with a convex reformulation. Section~\ref{s:5} explores the different risk-based optimization formulations for the thermal design of a cooling fin problem with non-convex limit state. Section~\ref{s:6} presents the concluding remarks.

\section{Reliability-based Design Optimization}\label{s:2} 
In this section, we review the RBDO formulation, which uses  PoF to quantify uncertainties. Let the quantity of interest of an engineering system be computed from the model $f: \mathcal{D}\times\Omega \mapsto \mathbb{R}$ as $f(\bd, Z)$, where the inputs to the system are the $n_d$ design variables $\bd\in\mathcal{D} \subseteq \mathbb{R}^{n_d}$ and the $n_z$ random variables $Z$ with the probability distribution $\pi$. The realizations of the random variables $Z$ are denoted by $\bz\in\Omega \subseteq \mathbb{R}^{n_z}$. The space of design variables is denoted by $\mathcal{D}$ and the space of random samples is denoted by $\Omega$. The failure of the system is described by a limit state function $g: \mathcal{D}\times\Omega \mapsto \mathbb{R}$ and a critical threshold $t \in \mathbb{R}$, where, without loss of generality, $g(\bd,\bz) > t$ defines failure of the system. For a system under uncertainty, $g(\bd,Z)$ is also a random variable given a particular design $\bd$. The limit state function in most engineering applications requires the solution of a system of equations (such as ordinary differential equations or partial differential equations).

\begin{figure}[!h]
	\centering
	\begin{tikzpicture}[scale=1,auto]
	\begin{axis}[
	xmin=0.3, xmax=3.5, ymin=0, ymax=1.0, height=3.5cm, width=7cm,
	black!60, very thick, domain=0:2.5, samples=100,
	axis lines=left,
	every axis y label/.style={anchor=south west,rotate=90,align=flush center,inner sep=0.3ex}, ylabel={Probability\\[-0.1cm] density},
	every axis x label/.style={at=(current axis.right of origin),anchor=north}, xlabel={$g(\bd,Z)$},  every tick/.style={black,	very thick}, tick label style={font=\selectfont,black}, xtick={2.2}, xticklabels={$t$}, ytick=\empty,
	enlargelimits=upper, clip=true, axis on top,
	smooth
	]
	\addplot [fill=red!60, draw=none, domain=2.2:3] {frechet(3,1.5,0)} \closedcycle;
	\addplot [very thick,black!80, domain=0.3:3] {frechet(3,1.5,0)};	
	
	\draw [decoration={markings,mark=at position 1 with {\arrow[scale=0.6]{*}}}, postaction={decorate},black!70]
	(axis cs:2.9,0.55) node [black, above, xshift=-0.3cm] 
	{$p_t=\mathbb{P} \left[g(\bd,Z)> t\right]$} -- (axis cs:2.5,0.07);
		
	\end{axis}		
	\end{tikzpicture}
	\caption{Illustration for PoF indicated by the area of the shaded region.}
	\label{fig:PoF}
\end{figure}
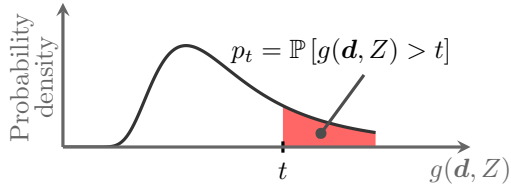

The most common RBDO formulation involves the use of a PoF constraint as
\begin{equation}
\label{e:RBDO}
\begin{split}
\min_{\bd\in\mathcal{D}}\quad & \mathbb{E}\left[f(\bd,Z)\right] \\
\text{subject to}\quad & \pt \left(g(\bd,Z)\right)\le 1-\alphaT,
\end{split}
\end{equation}
where $\alphaT \in[0,1]$ is the target reliability (i.e., $1-\alphaT$ is the target PoF) and the  PoF is defined via the limit state function $g$ and the failure threshold $t$ as $\pt\left(g(\bd,Z)\right):=\mathbb{P} \left[g(\bd,Z)> t\right]$. The RBDO problem~\eqref{e:RBDO} designs a system with optimal mean characteristics, in terms of $f(\bd,Z)$, such that it maintains a reliability of at least $\alpha$. Note, however, that PoF has no information about the magnitude of the failure event as it is merely a measure of the set $\{ g(\bd,Z)>t\}$, see Figure~\ref{fig:PoF}. 

For our upcoming discussion, it is helpful to point out that a constraint on the PoF is equivalent to a constraint on the $\alpha$-quantile. The $\alpha$-quantile, also known as the value-at-risk at level $\alpha$, is defined in terms of the inverse cumulative distribution function of the limit state function $F^{-1}_{g(\bd,Z)}$ as
\begin{equation}
\label{e:quant}
Q_\alpha\left[g(\bd,Z)\right] \coloneqq F_{g(\bd,Z)}^{-1}(\alpha).
\end{equation}
PoF and $Q_\alpha$ are natural counterparts that are measures of the tail of the distribution of $g(\bd,Z)$. When the largest $100(1-\alpha)\%$ outcomes are the ones of interest (i.e., failed cases), the quantile is a measure of minimum value within the set of these tail events. When one knows that outcomes larger than a given threshold $t$ are of interest, PoF provides a measure of the frequency of these ``large'' events. This equivalence of PoF and $Q_\alpha$ risk constraints is illustrated in Figure~\ref{fig:PoF_illustrate}. In the context of our optimization problem, using the same value of $t$ and $\alphaT$, \eqref{e:RBDO} can be written equivalently as
\begin{equation}
\label{e:RBDOvar}
\begin{split}
\min_{\bd\in\mathcal{D}}\quad & \mathbb{E}  \left[f(\bd,Z)\right] \\
\text{subject to}\quad & Q_{\alphaT}\left[g(\bd,Z)\right] \le t.
\end{split}
\end{equation}
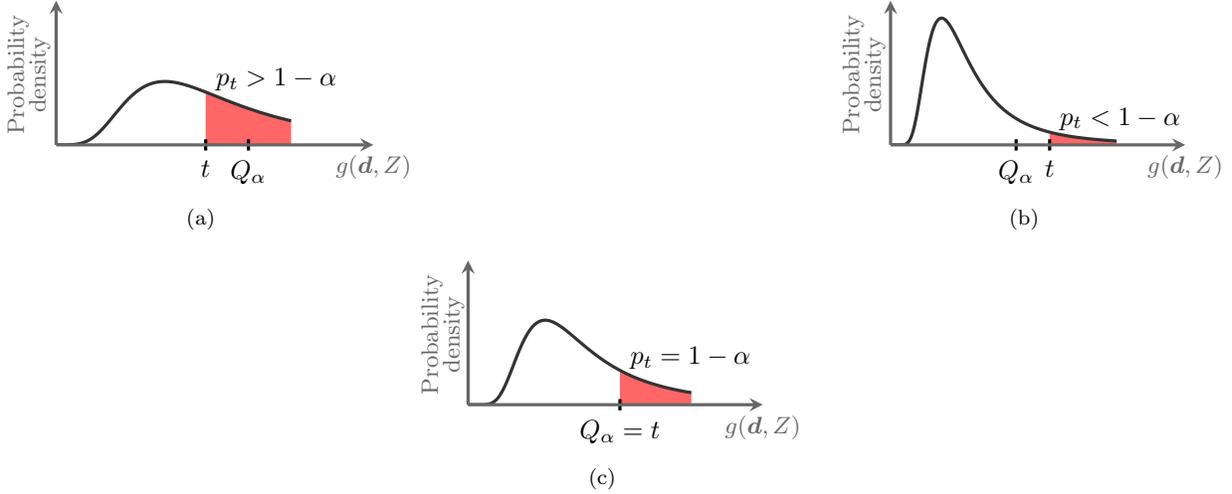
\begin{figure}[!htbp]
	\centering
	\subfigure[]{
		\begin{tikzpicture}[scale=1,auto]
		\begin{axis}[
		xmin=0.8, xmax=3.5, ymin=0, ymax=1.2, height=3.5cm, width=5.8cm,
		black!60, very thick, domain=0:2.5, samples=100,
		axis lines=left,
		every axis y label/.style={anchor=west,rotate=90,align=flush center,yshift=-0.35cm,font=\small}, ylabel={Probability\\[-0.1cm] density},
		every axis x label/.style={at=(current axis.right of origin),anchor=north,font=\small}, xlabel={$g(\bd,Z)$}, every tick/.style={black,	very thick}, tick label style={font=\selectfont,black}, xtick={2.2,2.6}, xticklabels={$t$,$Q_\alpha$}, ytick=\empty,
		enlargelimits=upper, clip=true, axis on top,
		smooth
		]
		\addplot [fill=red!60, draw=none, domain=2.2:3] {frechet(3,2,0)} \closedcycle;
		\addplot [very thick,black!80, domain=0.3:3] {frechet(3,2,0)};	
		
		\node[black,right] at (axis cs:2.2,0.6) {$p_t>1-\alpha$};
		\end{axis}		
		\end{tikzpicture}}\hfill
	\subfigure[]{
		\begin{tikzpicture}[scale=1,auto]
		\begin{axis}[
		xmin=0.3, xmax=3.5, ymin=0, ymax=1.2, height=3.5cm, width=5.5cm,
		black!60, very thick, domain=0:2.5, samples=100,
		axis lines=left,
		every axis y label/.style={anchor=west,rotate=90,align=flush center,yshift=0.35cm,font=\small}, ylabel={Probability\\[-0.1cm] density},
		every axis x label/.style={at=(current axis.right of origin),anchor=north,font=\small}, xlabel={$g(\bd,Z)$}, every tick/.style={black,	very thick}, tick label style={font=\selectfont,black}, xtick={1.8,2.2}, xticklabels={$Q_\alpha$,$t$}, ytick=\empty,
		enlargelimits=upper, clip=true, axis on top,
		smooth
		]
		\addplot [fill=red!60, draw=none, domain=2.2:3] {frechet(3,1,0)} \closedcycle;
		\addplot [very thick,black!80, domain=0.3:3] {frechet(3,1,0)};	
		
		\node[black,right] at (axis cs:2.2,0.24) {$p_t<1-\alpha$};
		\end{axis}		
		\end{tikzpicture}}\hfill
	\subfigure[]{
		\begin{tikzpicture}[scale=1,auto]
		\begin{axis}[
		xmin=0.5, xmax=3.5, ymin=0, ymax=1.2, height=3.5cm, width=5.5cm,
		black!60, very thick, domain=0:2.5, samples=100,
		axis lines=left,
		every axis y label/.style={anchor=west,rotate=90,align=flush center,yshift=0.1cm,font=\small}, ylabel={Probability\\[-0.1cm] density},
		every axis x label/.style={at=(current axis.right of origin),anchor=north,font=\small}, xlabel={$g(\bd,Z)$}, every tick/.style={black,	very thick}, tick label style={font=\selectfont,black}, xtick={2.2}, xticklabels={$Q_\alpha=t$}, ytick=\empty,
		enlargelimits=upper, clip=true, axis on top,
		smooth
		]
		\addplot [fill=red!60, draw=none, domain=2.2:3] {frechet(3,1.5,0)} \closedcycle;
		\addplot [very thick,black!80, domain=0.3:3] {frechet(3,1.5,0)};	
		
		\node[black,right] at (axis cs:2.2,0.43) {$p_t=1-\alpha$};
		\end{axis}		
		\end{tikzpicture}}
	\caption{Illustration of equivalence of PoF (shown by the shaded region) and $Q_\alpha$ showing that the two quantities converge at the constraint threshold when the reliability constraint is active.}
	\label{fig:PoF_illustrate}
\end{figure}

The most elementary method (although, inefficient) for estimating PoF is Monte Carlo~(MC) simulation when dealing with nonlinear limit state functions. The MC estimate of the PoF for a given design $\bd$ is
\begin{equation} \label{e:PFMC}
\hpt(g(\bd,Z)) = \frac{1}{m}\sum_{i=1}^{m}\mathbb{I}_{\failset(\bd)}(\bz_i),
\end{equation}
where $\bz_1,\dots,\bz_m$ are $m$ samples distributed according to $\pi$, $\failset(\bd) = \{ \bz \ | \ g(\bd,\bz) > t  \}$ is the failure set, and $\mathbb{I}_{\failset(\bd)}:\Omega \to \{0,1\}$ is the indicator function defined as
\begin{equation} \label{e:indicator}
\mathbb{I}_{\failset(\bd)}(\bz) = \left\{\begin{array}{ll}
1, & \text{if } \bz \in \failset(\bd) \\
0, & \text{else.}
\end{array}\right.
\end{equation}
The MC estimator is unbiased with the variance being $p_t(1-p_t)/m$.
The PoF estimation requires sampling from the tails of the distribution, which can often make MC estimators expensive.  A wealth of literature exists for methods that have been developed to deal with the computational complexity of PoF estimation and the RBDO problem. 
First, reliability index methods (e.g., FORM, SORM, etc.~\cite{sobieszczanski2015,du2001most}) geometrically approximate the limit state function to reduce the computational effort of PoF estimation. However, when the limit state function is nonlinear, the reliability index method could lead to inaccuracies in the estimate. 
Second, MC variance reduction techniques such as importance sampling~\cite{rubinstein2016simulation,melchers1989importance,mcbookArtOwen}, adaptive importance sampling~\cite{au1999new,dey1998ductile,de2005tutorial,papaioannou2015mcmc,depina2017coupling,kurtz2013cross}, and multifidelity approaches~\cite{dongbin2011,PKW17MFIS,PKW18MFCE} offer computational advantages. While the decay rate of the MC estimator cannot be improved upon, the variance of the MC estimator can be reduced, which offers computational advantages in that fewer (suitably chosen) MC samples are needed to obtain accurate PoF estimates. 
Third, adaptive data-driven surrogates for the limit state failure boundary identification can improve computational efficiency for the RBDO problem~\cite{bichon2013efficient,chaudhuri2019mfegra,moustapha2019surrogate}. 
Fourth, bi-fidelity RBDO methods~\cite{gano2006reliability,li2017vf} and recent multifidelity/multi-information-source methods for the PoF estimate~\cite{marques2018contour,KMPVW17FusionEstimators,chaudhuri2019mfegra} and the RBDO problem~\cite{chaudhuri2019reusing,chaudhuri2020information} have led to significant computational savings.

Although significant research has been devoted to PoF and RBDO, PoF as a risk measure does not factor in how catastrophic is the failure and thus, lacks resiliency. In other words, PoF neglects the magnitude of failure of the system and instead encodes a hard threshold via a binary function evaluation. We describe below this drawback of PoF.

\begin{remark}[Limitations of hard-thresholding] \label{r:hardthresh}
	To motivate the upcoming use of risk measures, we take a closer look at the limit state function $g$ and its use to characterize \textit{failure events}. In the standard setting, a failure event is characterized by a realization of $Z$ for some fixed design $\bd$ that leads to $g(\bd,\bz)>t$. However, this hard-threshold characterization of system failure potentially ignores important information quantified by the magnitude of $g(\bd,\bz)$ and PoF fails to promote resilience, i.e., no distinction between bad and very bad. Let us consider a structure with $g(\bd,\bz)$ being the load and the threshold $t$ being the allowable strength. There may be a large difference between the event $g(\bd,\bz)= t + .01kN$ and $g(\bd,\bz)= t + 100kN$, the latter characterizing a catastrophic system failure. This is not captured when considering system failure only as a binary decision with a hard threshold. Similarly, one could also consider events $g(\bd,\bz)= t - .01kN$ and $g(\bd,\bz)= t - 100kN$. A hard-threshold assessment deems both of these events as \textit{non-failure} events, even though $g(\bd,\bz)= t - .01kN$ is clearly a \textit{near-failure} event compared to $g(\bd,\bz)= t - 100kN$. A hard-threshold characterization of failure would potentially overlook these important \textit{near-failure} events and consider them as safe realizations of $g$. In reality, failure events do not usually occur using a hard-threshold rule. Even if they do, determination of the \textit{true} threshold will also involve uncertainty, blending statistical estimation, expert knowledge, and system models. Therefore, the choice of threshold should be involved in any discussion of measures of failure risk and we analyze later in Remark~\ref{r:conserve}, the advantage of the data-informed thresholding property of certain risk measures as compared to hard-thresholding. In addition, encoding magnitude of failure can help distinguish between designs with same PoF (see example 1 in Ref.~\cite{rockafellar2010buffered}). As we show in the next section, superquantile and bPoF do not have this deficiency. 
\end{remark} 

In the engineering community, PoF has been the preferred choice. Using PoF and RBDO offers some specific advantages starting with the simplicity of the risk measure and the natural intuition behind formulating the optimization problems, which is a major reason leading to the rich literature on this topic as noted before. Another advantage of PoF is the invariance to nonlinear reformulation for the limit state function. For example, let $z_1$ be a random load and $z_2$ be a random strength of a structure. Then the PoF would be the same regardless if the limit state function is defined as $z_1-z_2$ or $z_1/z_2-1$. Since $\alpha$-quantile leads to an equivalent formulation as PoF, both PoF and $\alpha$-quantile formulations have this invariance for continuous distributions.  However, there are several potential issues when using PoF as the risk measure for design optimization under uncertainty as noted below.

\begin{remark}[Optimization considerations] \label{r:discussionRBDO}
	While there are several advantages of using PoF and RBDO, there are several potential drawbacks. 
	First, PoF is not necessarily a convex function w.r.t.\ design variables $\bd$ even when the underlying limit state function is convex w.r.t.\ $\bd$. Thus, we cannot formulate a convex optimization problem even when underlying functions $f$ and $g$ are convex w.r.t.\ $\bd$. This is important because convexity guarantees convergence of standard and efficient algorithms to a globally optimal design under minimal assumptions since every local optimum is a global optimum in that case. 
	Second, the computation of PoF gradients can be ill-conditioned, so traditional gradient-based optimizers that require accurate gradient evaluations tend to face challenges. While PoF is differentiable for the specific case when $\bd$ only contains parameters of the distribution of $Z$, such as mean and standard deviation, PoF is in general not a differentiable function. Consequently, PoF gradients may not exist and when using approximate methods, such as finite difference, the accuracy of the PoF gradients could be poor. Some of these drawbacks can be addressed by using other methods for estimating the PoF gradients, but they have been developed under potentially restrictive assumptions~\cite{uryasev1995derivatives,royset2007extensions,tretiakov2000star}, which might not be easily verifiable for practical problems.
	Third, PoF can suffer from sensitivity to the failure threshold due to it being a discontinuous function w.r.t.\ threshold $t$. Since the choice of failure threshold could be uncertain, one would ideally prefer to have a measure of risk that is less sensitive to small changes in $t$. We further expand on this issue in Remark~\ref{r:contbPOF}.
\end{remark}

\section{Certifiable Risk-Based Design Optimization}\label{s:3}
Design optimization with a special class of risk measures can provide certifiable designs and algorithms.  We first present two notions of certifiability in risk-based optimization in Section~\ref*{s:3}.~\ref{s:CRiBDO}. We then discuss two specific risk measures, superquantile in Section~\ref*{s:3}.~\ref{s:SQopt} and bPoF in Section~\ref*{s:3}.~\ref{s:BPOF_opt}, that satisfy these notions of certifiability.

\subsection{Certifiability in risk-based design optimization}  \label{s:CRiBDO}
Risk in an engineering context can be quantified in several ways and the choice of risk measure, and its use as a cost or constraint, influences the design. 
We focus on a class of risk measures that can satisfy the following two \textit{certifiability conditions}:
\begin{enumerate}
	\item \textit{Data-informed conservativeness:}  Risk measures that take the magnitude of failure into account to decide the level of conservativeness required can certify the designs against near-failure and catastrophic failure events leading to increased resilience. The obtained designs can overcome the limitations of hard thresholding and are certifiably risk-averse against a continuous range of failure modes. In typical engineering problems, the limit state function distributions are not known and the information about the magnitude of failure is encoded through the generated data, thus making the conservativeness data-informed.
	\item \textit{Optimization convergence and efficiency:} Risk measures that preserve the convexity of underlying limit state functions (and/or cost functions) lead to convex risk-based optimization formulations. The resulting optimization problem is better behaved than a non-convex problem and can be solved more efficiently. Thus, one can find the design that is certifiably optimal in comparison with all alternate designs at reduced computational cost. In general, the risk measure preserves the convexity of the limit state function, such that the complexity of the optimization under uncertainty problem remains similar to the complexity of the deterministic optimization problem using the limit state function.
\end{enumerate}
We denote the risk-based design optimization formulations that use risk measures satisfying any of the two certifiability conditions as \textbf{C}ertifiable \textbf{Ri}sk-\textbf{B}ased \textbf{D}esign \textbf{O}ptimization (CRiBDO). 
Note that designs obtained through RBDO do not satisfy either of the above conditions since using PoF as the risk measure cannot guard against near-threshold or catastrophic failure events, see Remark~\ref{r:hardthresh}, and cannot certify the design to be a global optimum, see Remark~\ref{r:discussionRBDO}. Accounting for the magnitude of failure is critical to ensure appropriate conservativeness in CRiBDO designs and additionally, preservation of convexity is useful for optimizer efficiency and convergence guarantees to a globally optimal design. The optimization formulations satisfying both the conditions lead to \textit{strongly certifiable} risk-based designs. In general engineering applications, the convexity condition is difficult to satisfy but encapsulates an ideal situation, highlighting the importance of research in creating (piece-wise) convex approximations for physical problems.
In Sections~\ref*{s:3}.~\ref{s:SQopt} and ~\ref*{s:3}.~\ref{s:BPOF_opt}, we discuss the properties of two particular risk measures, superquantile and bPoF, that lead to certifiable risk-based designs and have the potential to be strongly certifiable when underlying functions are convex. Although we focus on these two particular risk measures in this work, other measures of risk could also be used to produce certifiable risk-based designs, see~\cite{artzner1999coherent,rockafellar2013fundamental,RTRockafellar_JORoyset_2015a}.

\subsection{Superquantile-based design optimization} \label{s:SQopt}
This section describes the concept of superquantiles and associated risk-averse optimization problem formulations. Superquantiles emphasize tail events, and from an engineering perspective it is important to manage such tail risks. 

\subsubsection{Risk measure: superquantile}\label{s:squant}
Intuitively, superquantiles can be understood as a tail expectation, or an average over a portion of worst-case outcomes. Given a fixed design $\bd$ and a distribution of potential outcomes $g(\bd,Z)$, the superquantile at level $\alpha\in[0,1]$ is the expected value of the largest $100(1-\alpha)\%$ realizations of $g(\bd,Z)$. In the literature, several other terms, such as CVaR and expected shortfall, have been used interchangeably with superquantile. We prefer the term superquantile because of its inherent connection with the long existing statistical quantity of quantiles and it being application agnostic. 

The definition of $\alpha$-superquantile is based on the $\alpha$-quantile $Q_\alpha\left[g(\bd,Z)\right]$ from Equation~\eqref{e:quant}. The $\alpha$-superquantile $\qbar_\alpha$ can be defined as
\begin{equation} \label{e:squant}
\qbar_\alpha\left[g(\bd,Z)\right] \coloneqq Q_\alpha\left[g(\bd,Z)\right] + \frac{1}{1-\alpha}\mathbb{E}  \left[\left[ g(\bd,Z)-Q_\alpha\left[g(\bd,Z)\right] \right]^+ \right],
\end{equation}
where $\bd$ is the given design and $[c]^+\coloneqq\max\{0,c\}$. The expectation in the second part of the right hand side of Equation~\eqref{e:squant} can be interpreted as the expectation of the tail of the distribution exceeding the $\alpha$-quantile. The $\alpha$-superquantile can be seen as the sum of the $\alpha$-quantile and a non-negative term and thus, $\qbar_\alpha \left[g(\bd,Z)\right]$ is a quantity higher (as indicated by ``super'') than $Q_\alpha\left[g(\bd,Z)\right]$. It follows from the definition that $\qbar_\alpha \left[g(\bd,Z)\right]\geq Q_\alpha\left[g(\bd,Z)\right]$. When the cumulative distribution of $g(\bd,Z)$ is continuous for any $\bd$, we can also view $\qbar_\alpha \left[g(\bd,Z)\right]$ as the conditional expectation of $g(\bd,Z)$ with the condition that $g(\bd,Z)$ is not less than $Q_\alpha\left[g(\bd,Z)\right]$, i.e., $\qbar_\alpha\left[g(\bd,Z)\right]=\mathbb{E}  \left[ g(\bd,Z)\ |\ g(\bd,Z)\ge Q_\alpha\left[g(\bd,Z) \right] \right ]$~\cite{rockafellar2000optimization}. We also note that by definition~\cite{rockafellar2013superquantiles}
\begin{equation}
\label{e:squantLim}
\begin{split}
&\text{for } \alpha=0,\ \qbar_0\left[g(\bd,Z)\right]=\mathbb{E}  \left[g(\bd,Z)\right], \text{ and } \\ 
&\text{for } \alpha=1,\ \qbar_1\left[g(\bd,Z)\right]=\operatornamewithlimits{ess\,sup} g(\bd,Z),
\end{split}
\end{equation}
where $\operatornamewithlimits{ess\,sup} g(\bd,Z)$ is the essential supremum, i.e., the lowest value that $g(\bd,Z)$ doesn't exceed with probability 1. 

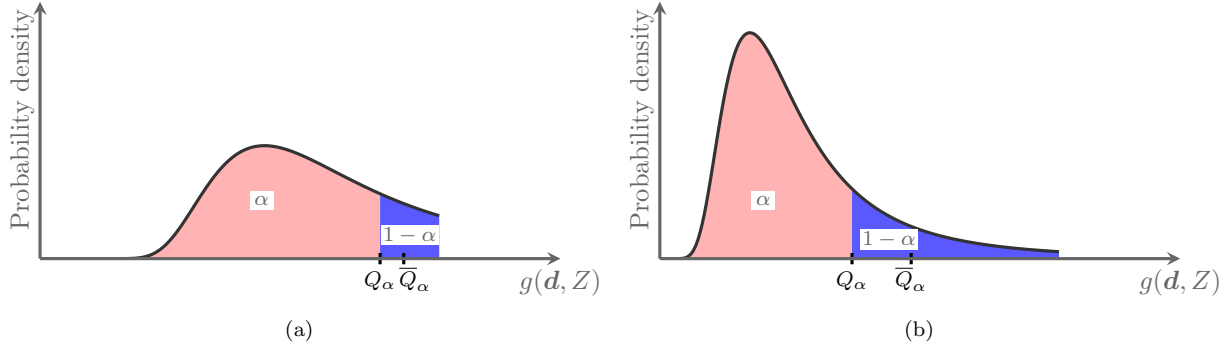
\begin{figure}[!htbp]
	\centering
	\subfigure[]{
		\begin{tikzpicture}[scale=1,auto]
		\begin{axis}[
		xmin=0.3, xmax=3.5, ymin=0, ymax=1.2, height=5cm, width=8.5cm,
		black!60, very thick, domain=0:2.5, samples=100,
		axis lines=left,
		every axis y label/.style={anchor=south west,rotate=90,align=flush center,inner sep=0.4ex,xshift=0.3cm}, ylabel={Probability density},
		every axis x label/.style={at=(current axis.right of origin),anchor=north}, xlabel={$g(\bd,Z)$}, every tick/.style={black,	very thick}, tick label style={font=\footnotesize,black,text height=1ex}, xtick={2.6,2.76}, xticklabels={$Q_{\alpha}$,\hspace*{2.5ex}$\qbar_{\alpha}$}, ytick=\empty, enlargelimits=upper, clip=true, axis on top, smooth
		]
		\addplot [fill=red!30, draw=none, domain=0.3:2.6] {frechet(3,2,0)} \closedcycle;
		\addplot [fill=blue!65, draw=none, domain=2.6:3] {frechet(3,2,0)} \closedcycle;
		\addplot [very thick,black!80, domain=0.3:3] {frechet(3,2,0)};	
		
		\node at (axis cs:1.8,0.3) [fill=white,font=\footnotesize,inner sep=0.4ex] {$\alpha$};	
		\node at (axis cs:2.8,0.12) [fill=white,font=\footnotesize,inner sep=0.18ex] {$1-\alpha$};
		\end{axis}		
		\end{tikzpicture}}\hfill
	\subfigure[]{
		\begin{tikzpicture}[scale=1,auto]
		\begin{axis}[
		xmin=0.3, xmax=3.5, ymin=0, ymax=1.2, height=5cm, width=8.5cm,
		black!60, very thick, domain=0:2.5, samples=100,
		axis lines=left,
		every axis y label/.style={anchor=south west,rotate=90,align=flush center,inner sep=0.4ex,xshift=0.3cm}, ylabel={Probability density},
		every axis x label/.style={at=(current axis.right of origin),anchor=north}, xlabel={$g(\bd,Z)$}, every tick/.style={black,	very thick}, tick label style={font=\footnotesize,black,text height=1ex}, xtick={1.6,2.0}, xticklabels={$Q_\alpha$,$\qbar_{\alpha}$}, ytick=\empty,
		enlargelimits=upper, clip=true, axis on top, smooth
		]
		\addplot [fill=red!30, draw=none, domain=0.3:1.6] {frechet(3,1,0)} \closedcycle;
		\addplot [fill=blue!65, draw=none, domain=1.6:3] {frechet(3,1,0)} \closedcycle;
		\addplot [very thick,black!80, domain=0.3:3] {frechet(3,1,0)};
		
		\node at (axis cs:1,0.3) [fill=white,font=\footnotesize,inner sep=0.4ex] {$\alpha$};	
		\node at (axis cs:1.85,0.1) [fill=white,font=\footnotesize,inner sep=0.2ex] {$1-\alpha$};
		\end{axis}		
		\end{tikzpicture}}
	\caption{Illustration for $\qbar_\alpha$ on two generic distributions: expectation of the worst-case $1-\alpha$ outcomes shown in blue is $\qbar_\alpha \left[g(\bd,Z)\right]$.}
	\label{fig:CVaR_illustrate}
\end{figure}

Figure~\ref{fig:CVaR_illustrate} illustrates the $\qbar_\alpha$ risk measure for two differently shaped, generic distributions of the limit state function. The figure shows that the magnitude of $\qbar_{\alpha}-Q_\alpha$ (or the induced conservativeness) changes with the underlying distribution. Algorithm~\ref{algo:VaR_and_CVaR_from_samples} describes standard MC sampling for approximating $\qbar_\alpha$. The second term on the right hand side in Equation~\eqref{e:CVaR_discrete} is a MC estimate of the expectation in Equation~\eqref{e:squant}.
%
\begin{algorithm}[!htb]
	\caption{Sampling-based estimation of $Q_\alpha$ and $\qbar_\alpha$.}
	\label{algo:VaR_and_CVaR_from_samples}
	\begin{algorithmic}[1]
		\Require $m$ i.i.d. samples $\bz_1,\dots,\bz_m$ of random variable $Z$, design variable $\bd$, risk level $\alpha \in (0,1)$, limit state function  $g(\bd,Z)$.
		\Ensure Sample approximations $\widehat{Q}_\alpha \left[g(\bd,Z)\right]$, $\widehat{\qbar}_\alpha \left[g(\bd,Z)\right]$.
		
		\State Evaluate limit state function at the samples to get  $g(\bd,\bz_1),\dots,g(\bd,\bz_m)$.
		
		\State Sort values of limit state function in descending order and relabel the samples so that
		\[
		g(\bd,\bz_1) >  g(\bd,\bz_2) > \cdots > g(\bd,\bz_m).
		\]
		
		\State Find the index $k_\alpha = \lceil m(1-\alpha)\rceil$ to estimate $\widehat{Q}_\alpha\left[g(\bd,Z)\right] \gets  g(\bd,\bz_{k_\alpha})$.
		
		\State Estimate
		\begin{equation} \label{e:CVaR_discrete}
		\widehat{\qbar}_\alpha \left[g(\bd,Z)\right]  = \widehat{Q}_\alpha\left[g(\bd,Z)\right] + \frac{1}{m(1-\alpha)} \sum_{j=1}^{m} \left[g(\bd,\bz_j)-\widehat{Q}_\alpha\left[g(\bd,Z)\right]\right]^+.
		\end{equation}
	\end{algorithmic}
\end{algorithm}

\subsubsection{Optimization problem: superquantiles as constraint}
As noted before, the PoF constraint of the RBDO problem in \eqref{e:RBDO} can be viewed as a $Q_\alpha$ constraint (as seen in \eqref{e:RBDOvar}). The PoF constraint (and thus the $Q_\alpha$ constraint) does not consider the magnitude of the failure events, but only whether they are larger than the failure threshold. This could be a potential drawback for engineering applications. On the other hand, a $\qbar_\alpha$ constraint considers the magnitude of the failure events by specifically constraining the expected value of the largest $100(1-\alpha)\%$ realizations of $g(\bd,Z)$. Additionally, depending upon the actual construction of $g(\bd,\bz)$ and the accuracy of the sampling procedure, the $\qbar_\alpha$ constraint may have numerical advantages over the $Q_\alpha$ constraint when it comes to optimization as discussed later. In particular, we have the optimization problem formulation
\begin{equation}
\label{e:CVaRopt1}
\begin{split}
\min_{\bd\in\mathcal{D}}\quad & \mathbb{E}  \left[f(\bd,Z)\right] \\
\text{subject to}\quad & \qbar_{\alphaT} \left[g(\bd,Z)\right]\le t,
\end{split}
\end{equation}
where $\alphaT$ is the desired reliability level given the limit state failure threshold $t$. The $\qbar_\alpha$-based formulation typically leads to a more conservative design than when PoF is used. This can be observed by noting that $\qbar_{\alphaT}\left[g(\bd,Z)\right] \le t \implies Q_{\alphaT}\left[g(\bd,Z)\right]\le t \iff  \pt(g(\bd,Z))\le 1-\alphaT$. Therefore, if the design satisfies the $\qbar_{\alphaT}$ constraint, then the design will also satisfy the related PoF constraint. Additionally, since the $\qbar_{\alphaT}$ constraint ensures that the average of the $(1-\alphaT)$ tail is no larger than $t$, it is likely that the probability of exceeding $t$ (PoF) is strictly smaller than $1-\alphaT$ and is thus a conservative design for target reliability of $\alphaT$. Intuitively, this conservatism comes from the fact that $\qbar_{\alphaT}$ considers the magnitude of the worst failure events. 

The formulation with $\qbar_{\alphaT}$ as the constraint is useful when the designer is unsure about the failure boundary location for the problem but requires a certain level of reliability from the design. For example, consider the case where the failure is defined as maximum stress of a structure not exceeding a certain value. However, the designers cannot agree on the cut-off value for stress but can agree on the desired level of reliability they want. One can use this formulation to design a structure with a given reliability ($1-\alphaT$) while constraining a conservative estimate of the cut-off value ($\qbar_{\alphaT}$) on the stress.

\begin{remark}[Convexity in $\qbar_\alpha$-based optimization] \label{r:convexCVaR}
	It can be shown that $\qbar_\alpha$ can be written in the form of an optimization problem~\cite{rockafellar2000optimization} as 
	\begin{equation} \label{e:squantOpt}
	\qbar_\alpha \left[g(\bd,Z)\right] = \min_{\gamma\in\mathbb{R}}  \gamma+ \frac{1}{1-\alpha}\mathbb{E}  \left[ \left[g(\bd,Z) -\gamma \right]^+\right],
	\end{equation}
	where $\bd$ is the given design, $\gamma$ is an auxiliary variable, and $[c]^+\coloneqq\max\{0,c\}$. At the optimum, $\gamma^* = Q_\alpha\left[g(\bd,Z)\right]$.
	Using Equation~\eqref{e:squantOpt}, the formulation~\eqref{e:CVaRopt1} can be reduced to an optimization problem involving only expectations as given by
	\begin{equation}
	\label{e:squantopt1_gammma}
	\begin{split}
	\min_{\gamma\in\mathbb{R},\, \bd\in\mathcal{D}}\quad & \mathbb{E}  \left[f(\bd,Z)\right] \\
	\emph{subject to} \quad &  \gamma+ \frac{1}{1-\alphaTrem}\mathbb{E}  \left[ \left[ g(\bd,Z) -\gamma \right]^+\right] \le t.
	\end{split}
	\end{equation}
	The formulation \eqref{e:squantopt1_gammma} is a convex optimization problem when $g(\bd,Z)$ and $f(\bd,Z)$ are convex in $\bd$ since $[\cdot]^+$ is a convex function and preserves the convexity of the limit state function.
	Another advantage of \eqref{e:squantopt1_gammma}, as outlined in Ref.~\cite{rockafellar2000optimization}, is that the nonlinear part of the constraint, $\mathbb{E}  \left[ \left[g(\bd,Z) -\gamma \right]^+\right]$, can be reformulated as a set of convex (linear) constraints if $g(\bd,Z)$ is convex (linear) in $\bd$ and has a discrete (or empirical) distribution with the distribution of $Z$ being independent of $\bd$~\footnote{In the case where $\pi$ depends upon $\bd$, one can perform optimization by using sampling-based estimators for the gradient of $\qbar_\alpha$~\cite{lan2016algorithms,tamar2015optimizing}.}.  
	Specifically, consider a MC estimate where $\bz_i, i=1,\dots,m$ are $m$ samples from probability distribution $\pi$. Then, using auxiliary variables $b_i, i=1,\dots,m$ to define $\bb=\{b_1,\dots,b_m\}$, we can reformulate \eqref{e:squantopt1_gammma} as
	\begin{equation}
	\label{e:squantopt1_gammma_linear}
	\begin{split}
	\min_{\gamma\in\mathbb{R},\, \bb\in\mathbb{R}^m,\, \bd\in\mathcal{D}}\quad & \mathbb{E}  \left[f(\bd,Z)\right] \\
	\emph{subject to}\quad &  \gamma+ \frac{1}{m(1-\alphaTrem)} \sum_{i=1}^m b_i \le t, \\
	& g(\bd,z_i)-\gamma \le b_i, i=1,\dots,m, \\
	& b_i \ge 0, i=1,\dots,m.
	\end{split}
	\end{equation}
	The formulation \eqref{e:squantopt1_gammma_linear} is a linear program when $g(\bd,Z)$ and $f(\bd,Z)$ are linear in $\bd$.
\end{remark}

As noted in Remark~\ref{r:convexCVaR}, the formulations in \eqref{e:squantopt1_gammma} and \eqref{e:squantopt1_gammma_linear} are convex (or linear) only when the underlying functions $g(\bd,Z)$ and $f(\bd,Z)$ are convex (or linear) in $\bd$. However, the advantages and possibility of such formulations indicates that one can achieve significant gains by investing in convex (or linear) approximations for the underlying functions.

\subsubsection{Optimization problem: superquantiles as objective}
The $\alpha$-superquantile $\qbar_\alpha$ naturally arises as a replacement for $Q_\alpha$ in the constraint, but it can also be used as the objective function in the optimization problem formulation. For example, in PDE-constrained optimization, superquantiles have been used in the objective function~\cite{DPKouri_TMSurowiec_2016a,ZZou_DPKouri_WAquino_2018a}. The optimization formulation is
\begin{equation}
\label{e:CVaRopt2}
\begin{split}
\min_{\bd\in\mathcal{D}}\quad &  \qbar_{\alphaT} \left[g(\bd,Z)\right]\\
\text{subject to}\quad & \qbar_{\beta_{\scriptscriptstyle \text{T}}} \left[f(\bd,Z)\right] \le C_{\scriptscriptstyle \text{T}},
\end{split}
\end{equation}
where $\alphaT$ and $\beta_{\scriptscriptstyle \text{T}}$ are the desired risk levels for $g$ and $f$ respectively, and $C_{\scriptscriptstyle \text{T}}$ is a threshold on the quantity of interest $f$. This is a useful formulation when it is easier to define a threshold on the quantity of interest than deciding a risk level for the limit state function. For example, if the quantity of interest is the cost of manufacturing a rocket engine, one can specify a budget constraint and use the above formulation. The solution of this optimization formulation would result in the safest rocket engine design such that the expected budget does not exceed the given budget.

\subsubsection{Discussion on superquantile-based optimization} \label{s:cvarDisc}
From an optimization perspective, an important feature of $\qbar_\alpha$ is that it preserves convexity of the function it is applied to, i.e., the limit state function or cost function. $\qbar_\alpha$-based formulations can lead to well-behaved convex optimization problems that allows one to provide convergence guarantees as described in Remark~\ref{r:convexCVaR}. The reformulation offers a major advantage, since an optimization algorithm can work directly on the limit state function without passing through an indicator function. This preserves the convexity and other mathematical properties of the limit state function. $\qbar_\alpha$ also takes the magnitude of failure into account, which makes it more informative and resilient compared to PoF and builds in data-informed conservativeness.

As noted in \cite{yamai2002comparative}, $\qbar_\alpha$ estimators are less stable than estimators of $Q_\alpha$ since rare, large magnitude tail samples can have large effect on the sample estimate. This is more prevalent when the distribution of the random quantity is fat-tailed. Thus, there is a need for more research to develop efficient algorithms for $\qbar_\alpha$ estimation.
Despite offering convexity, a drawback of $\qbar_\alpha$ is that it is non-smooth, and a direct $\qbar_\alpha$-based optimization would require either non-smooth optimization methods, for example variable-metric algorithms \cite{uryas1991new}, or gradient-free methods. Note that smoothed approximations exist~\cite{DPKouri_TMSurowiec_2016a,basova2011computational}, which significantly improve optimization performance. In addition, the formulation  \eqref{e:squantopt1_gammma_linear} offers a smooth alternative.

As noted in Remark~\ref{r:convexCVaR}, $\qbar_\alpha$-based formulations can be further reduced to a linear program. The formulation in \eqref{e:squantopt1_gammma_linear} increases the dimensionality of the optimization problem from $n_d+1$ to $n_d+m+1$, where $m$ is the number of MC samples, which poses an issue when the number of MC samples is large. However, formulation \eqref{e:squantopt1_gammma_linear} has mostly linear constraints and can also be completely converted into a linear program by using a linear approximation for $g(\bd,z_i)$ (following similar ideas as reliability index methods described in Section~\ref{s:introduction}). There are extremely efficient methods for finding solutions to linear programs even for high-dimensional problems.

\subsection{bPoF-based design optimization} \label{s:BPOF_opt} 
Buffered probability of failure was first introduced by Rockafellar and Royset~\cite{rockafellar2010buffered} as an alternative to PoF. This section describes bPoF and the associated optimization problem formulations. When used as constraints, bPoF and superquantile lead to equivalent optimization formulations but bPoF provides an alternative interpretation of the $\qbar_\alpha$ constraint that is, arguably, more natural for applications dealing with constraints in terms of failure probability instead of constraints involving quantiles. When considered as an objective function, bPoF and superquantile lead to different optimal design solutions.

\subsubsection{Risk measure: bPoF}
The bPoF is an alternate measure of reliability which adds a buffer to the traditional PoF. The definition of bPoF at a given design $\bd$ is based on the superquantile as given by
\begin{equation}
\label{e:bPoF_alt_1}
\pbuf \left(g(\bd,Z)\right)\coloneqq \left\{\begin{array}{lll}
\left\{ 1-\alpha \; | \; \qbar_{\alpha}\left[g(\bd,Z)\right]=t \right\}, &\text{if} \ \qbar_0\left[g(\bd,Z)\right] < t < \qbar_1\left[g(\bd,Z)\right] \\
0, &\text{if} \ t \geq \qbar_1\left[g(\bd,Z)\right] \\
1, &\text{otherwise}.
\end{array}\right. 
\end{equation}
The domains of the threshold $t$ in Equation~\eqref{e:bPoF_alt_1} can interpreted in more intuitive terms using Equation~\eqref{e:squantLim} for $\qbar_0\left[g(\bd,Z)\right]$ and $\qbar_1\left[g(\bd,Z)\right]$. 
The relationship between superquantiles and bPoF in the first condition in Equation~\eqref{e:bPoF_alt_1} can also be viewed in the same way as that connecting $\alpha$-quantile and PoF by recalling that
\begin{equation}
\label{e:qbarbPoFequiv}
Q_{\alpha}\left[g(\bd,Z)\right] \leq t \iff  \pt(g(\bd,Z)) \leq 1-\alpha \text{ and here, }\ \qbar_{\alpha}\left[g(\bd,Z)\right] \leq t \iff \pbuf \left(g(\bd,Z)\right) \leq 1 - \alpha.
\end{equation} 

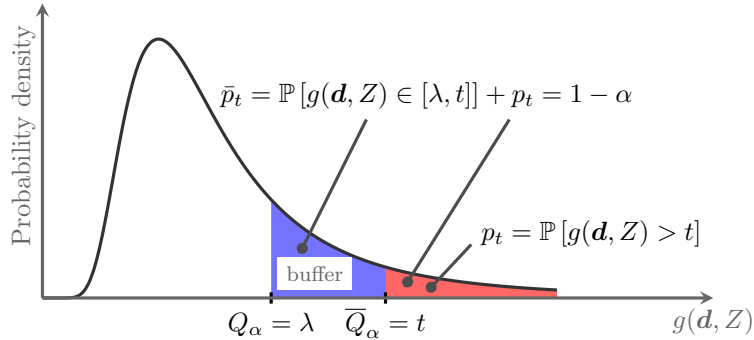
\begin{figure}[!htb]
	\centering
	\begin{tikzpicture}[scale=1,auto]
	
	\begin{axis}[
	xmin=0.3, xmax=3.5, ymin=0, ymax=1.2, height=5.5cm, width=10.5cm,
	black!60, very thick, domain=0:2.5, samples=100,
	axis lines=left,
	every axis y label/.style={anchor=south west,rotate=90,align=flush center,inner sep=0.5ex,xshift=0.6cm}, ylabel={Probability density},
	every axis x label/.style={at=(current axis.right of origin),anchor=north}, xlabel={$g(\bd,Z)$}, every tick/.style={black,	very thick}, tick label style={font=\selectfont,black,text height=1.5ex}, xtick={1.5,2.1}, xticklabels={$Q_{\alpha}=\lambda$,$\qbar_{\alpha}=t$}, ytick=\empty,
	enlargelimits=upper, clip=false, axis on top,
	smooth
	]
	\addplot [fill=blue!55, draw=none, domain=1.5:3] {frechet(3,1,0)} \closedcycle;
	\addplot [fill=red!60, draw=none, domain=2.1:3] {frechet(3,1,0)} \closedcycle;
	\addplot [very thick,black!80, domain=0.3:3] {frechet(3,1,0)};
	
	\node at (axis cs:1.73,0.11) [fill=white,font=\footnotesize] {buffer};	
	
	\draw [decoration={markings,mark=at position 1 with {\arrow[scale=0.6]{*}}}, postaction={decorate},black!70]
	(axis cs:2.55,0.24) node [black, right, yshift=0.15cm] 
	{$p_t=\mathbb{P} \left[g(\bd,Z)>t\right]$} -- (axis cs:2.32,0.03);
	
	
	\draw [decoration={markings,mark=at position 1 with {\arrow[scale=0.6]{*}}}, postaction={decorate},black!70]
	(axis cs:2,0.8) node [black, above, xshift=0.8cm] 
	{$\bar{p}_t = \mathbb{P} \left[g(\bd,Z)\in[\lambda,t]\right] + p_t =1-\alpha$} -- (axis cs:1.65,0.2);
	\draw [decoration={markings,mark=at position 1 with {\arrow[scale=0.6]{*}}}, postaction={decorate},black!70](axis cs:2.75,0.8) -- (axis cs:2.2,0.05);
	\end{axis}
	
	\end{tikzpicture}
	\caption{Illustration for bPoF: for a given threshold $t$, PoF equals the area in \textit{red} while bPoF equals the combined area in \textit{red and blue}.}
	\label{fig:bPOF_illustrate}
\end{figure}

To make the concept of \textit{buffer} concrete, we further analyze the case in the first condition in Equation~\eqref{e:bPoF_alt_1} when $t \in \left(\qbar_0\left[g(\bd,Z)\right],\qbar_1\left[g(\bd,Z)\right]\right)$ and $g(\bd,Z)$ is a continuous random variable, which leads to $\pbuf \left(g(\bd,Z)\right)=\left\{ 1-\alpha \; | \; \qbar_{\alpha}\left[g(\bd,Z)\right]=t \right\}$. Using the definition of quantiles from Equation~\eqref{e:quant} and its connection with superquantiles (see Equation~\eqref{e:squant} and Figure~\ref{fig:CVaR_illustrate}), we can see that $1-\alpha=\mathbb{P}\left[g(\bd,Z)\ge Q_{\alpha}[g(\bd,Z)]\right]$. This leads to another definition of bPoF in terms of probability of exceeding a quantile given the condition on $\alpha$ as
\begin{equation} \label{e:bPoF}
\pbuf \left(g(\bd,Z)\right)=\mathbb{P}  \left[g(\bd,Z)\ge Q_{\alpha}\left[g(\bd,Z)\right]\right]=1-\alpha, \quad  \text{where $\alpha$ is such that} \quad  \qbar_{\alpha}\left[g(\bd,Z)\right]=t.
\end{equation}
We know that superquantiles are conservative as compared to quantiles (Section~\ref*{s:3}.~\ref*{s:SQopt}.~\ref{s:squant}), which leads to $Q_{\alpha}\le t$ since $\qbar_{\alpha}=t$. Thus, Equation~\eqref{e:bPoF} can be split as a sum of PoF and the probability of \textit{near-failure} as
\begin{equation}
\label{e:bPoFsum}
\pbuf \left(g(\bd,Z)\right) = \mathbb{P}\left[g(\bd,Z) > t \right] + \mathbb{P}\left[ g(\bd,Z) \in \left[Q_{\alpha}\left[g(\bd,Z)\right],t\right] \right] = \pt\left(g(\bd,Z)\right) + \mathbb{P}\left[ g(\bd,Z) \in \left[\lambda,t\right] \right], 
\end{equation}
where $\lambda=Q_{\alpha}\left[g(\bd,Z)\right]$. The value of $\lambda$ is affected by the condition on $\alpha$ through superquantiles (see Equation~\eqref{e:bPoF}) and takes into account the frequency and magnitude of failure. Thus, the near-failure region $\left[\lambda,t\right]$ is determined by the frequency and magnitude of tail events around $t$ and can be intuitively seen as the \textit{buffer} on top of the PoF. An illustration of the bPoF risk measure is shown in Figure~\ref{fig:bPOF_illustrate}. Algorithm~\ref{algo:bPoF_from_samples} describes standard MC sampling for estimating bPoF. Note that all the quantities discussed in this work (PoF, superquantiles, and bPoF) can be viewed as expectations. Estimating them via Monte Carlo simulation therefore yields estimates whose error decreases with the rate $ 1/\sqrt{\text{number of samples}} $. All the estimates suffer from an increasing constant associated with the estimator variance as one moves further out in the tail, i.e., larger threshold or larger $\alpha$. The computational effort can be reduced for any of the risk measures by using Monte Carlo variance reduction strategies.
\begin{algorithm}[!htb]
	\caption{Sampling-based estimation of bPoF.}
	\label{algo:bPoF_from_samples}
	\begin{algorithmic}[1] 
		\Require $m$ i.i.d. samples $\bz_1,\dots,\bz_m$ of random variable $Z$, design variable $\bd$, failure threshold $t$, and limit state function  $g(\bd,Z)$.
		\Ensure Sample approximation $\widehat{\pbuf} \left(g(\bd,Z)\right)$.
		
		\State Evaluate limit state function at the samples to get $g(\bd,\bz_1),\dots,g(\bd,\bz_m)$.
		
		\State Sort values of limit state function in descending order and relabel the samples so that
		\[
		g(\bd,\bz_1) >  g(\bd,\bz_2) > \ldots > g(\bd,\bz_m).
		\]
		\State $c=g(\bd,\bz_1)$	\Comment{Initialize superquantile estimate}
		\State $k=1$
		\While{$c\ge t$}	\Comment{Check if superquantile estimate equals threshold}
		\State $k \gets k+1$
		\State $c=\frac{1}{k} \sum_{i=1}^k g(\bd,\bz_k)$	\Comment{Update superquantile estimate}
		\EndWhile
		\State Estimate bPoF as $\widehat{\pbuf} \left(g(\bd,Z)\right) \approx \frac{k-1}{m}  $ \Comment{Estimate bPoF as $1-\alpha$ when $c\approx t$}
	\end{algorithmic}
\end{algorithm} 

In general, we can see that for any design $\bd$,
\begin{equation} \label{e:bPOFconserve}
\pbuf \left(g(\bd,Z)\right)\ge \pt\left(g(\bd,Z)\right).
\end{equation}
Through Equation~\eqref{e:bPoFsum}, we can see that the conservatism of bPoF comes from the data-dependent mechanism that selects the conservative threshold $\lambda \le t$, which acts to establish a \textit{buffer} zone.  
If realizations of $g(\bd,Z)$ beyond $t$ are very large (potentially catastrophic failures), $\lambda$ will need to be smaller (making bPoF bigger) to drive the expectation beyond $\lambda$ to $t$. Thus, the larger bPoF serves to account for not only the frequency of failure events, but also their magnitude. The bPoF also accounts for the frequency of \textit{near-failure} events that have magnitude below, but very close to $t$. If there are a large number of near-failure events, bPoF will take this into account, since it will be included in the $\lambda$-tail which must have average equal to $t$.
Thus, the bPoF is a conservative estimate of the PoF for any design $\bd$ and carries more information about failure than PoF since it takes into consideration the magnitude of failure. 
It has been shown that for exponential distribution of the limit state function, the bPoF is $e\approx2.718$ times the PoF~\cite{pertaia2021new}. However, the degree of conservativeness of bPoF w.r.t. PoF is dependent on the distribution of $g(\bd,Z)$, which is typically not known in closed-form in engineering applications.

\begin{remark}[Continuity of bPoF w.r.t.\ threshold] \label{r:contbPOF}
	In practice, thresholds are sometimes set by regulatory commissions, informed by industry standards; see chapter 18 in Ref.~\cite{ditlevsen1996structural} for a discussion on code calibration.
	As discussed before, the data-informed conservativeness of bPoF reduces the adverse effects of poorly chosen thresholds by building a buffer around the threshold $t$. Another issue with poorly set thresholds is that the values could change as one learns more about the system. In such cases, continuity of the risk measure w.r.t.\ the threshold becomes important. bPoF is continuous w.r.t.\ the threshold but PoF is not. Consequently, if an engineer makes small changes to the threshold $t$, then it can have significant effects on the resulting design when PoF is used in the optimization formulation. On the other hand, small changes in $t$ will only have small effect on the bPoF-based optimal design due to bPoF being continuous w.r.t.\ $t$.\\
	The following example illustrates the continuity of bPoF vs PoF w.r.t.\ the threshold. Let $X$ be a random variable with finite distribution probability mass function given by
	\begin{equation*}
	\mathbb{P}(X=x) = \left\{\begin{array}{lll}
	0.8, &\text{if } x=-1 \\
	0.1, &\text{if } x=0 \\
	0.1, &\text{if } x=1,
	\end{array}\right.
	\end{equation*}
	which is visualized in Figure~\ref{fig:ExcontThresh}(a). For this simple distribution, one can derive the PoF and bPoF analytically for any given threshold $t$. The PoF values for different values of $t$ are
	\begin{equation*}
	p_t = \left\{\begin{array}{llll}
	1, &\text{if } t<-1 \\
	0.2, &\text{if } t\in[-1,0) \\
	0.1, &\text{if } t\in[0,1)\\
	0, &\text{if } t\ge 1,
	\end{array}\right.
	\end{equation*}
	which is clearly not continuous in $t$. The bPoF values for different values of $t$ are
	\begin{equation*}
	\pbuf = \left\{\begin{array}{llll}
	1, &\text{if } t < -0.7 \\
	0.3/(t+1), &\text{if } t\in[-0.7,0.5) \\
	0.1/t, &\text{if } t\in[0.5,1)\\
	0, &\text{if } t \ge 1,
	\end{array}\right.
	\end{equation*}
	which is continuous in $t$ on the interval $\left(-\infty,\qbar_1\right)=\left(-\infty,1 \right)$. The PoF and bPoF values as a function of the threshold $t$ are plotted in Figure~\ref{fig:ExcontThresh}(b) showing the continuity of bPoF in $t$. In a similar way, superquantiles $\qbar_{\alpha}$ are continuous in $\alpha$ but quantiles $Q_\alpha$ are not.
	\begin{figure}[!htb]
		\centering
		\subfigure[]{\includegraphics[width=0.5\textwidth,page=1]{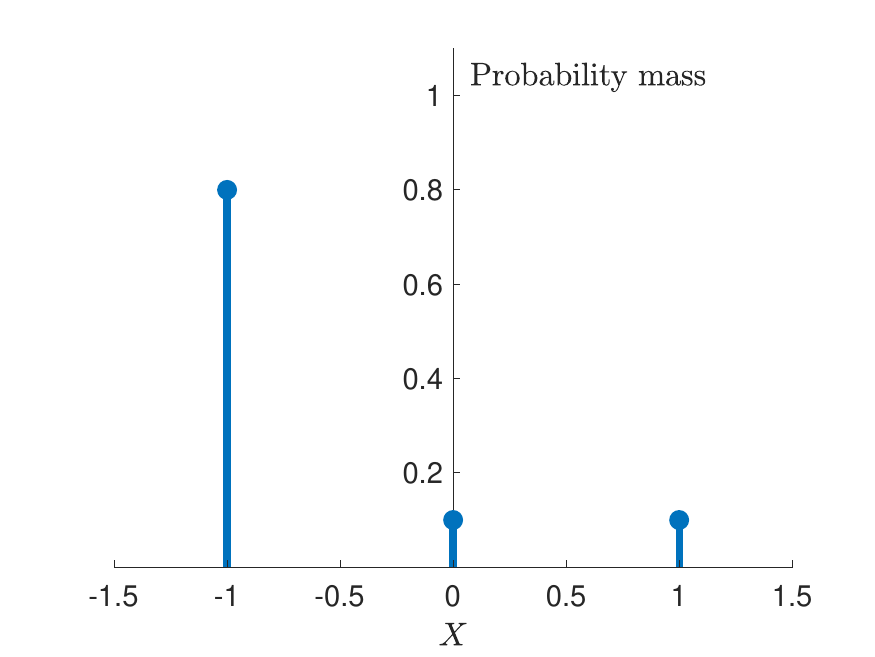}}\hfill
		\subfigure[]{\includegraphics[width=0.5\textwidth,page=2]{PoFvsbPoF_contEx.pdf}}
		\caption{Illustrating (a) the probability mass function of $X$ and (b) the continuity of bPoF w.r.t.\ changing threshold values as compared to discontinuous nature of PoF.}
		\label{fig:ExcontThresh}
	\end{figure}
\end{remark}

\subsubsection{Optimization problem: bPoF as constraint}
One of the advantages of bPoF, which provides data-informed conservativeness, is the intuitive relatability to the widely used PoF. This helps in easy transition from PoF-based formulations to bPoF-based formulations. 
Consider the optimization problem \eqref{e:RBDO} with the PoF constraint replaced by the bPoF constraint,
\begin{equation} 
\label{e:bPOFopt1}
\begin{split}
\min_{\bd\in\mathcal{D}}\quad & \mathbb{E}  \left[f(\bd,Z)\right] \\
\text{subject to}\quad & \pbuf \left(g(\bd,Z)\right)\le 1-\alphaT.
\end{split}
\end{equation}
Just as a PoF constraint is equivalent to a $Q_\alpha$ constraint, it can be shown that the bPoF constraint formulation described above is equivalent to a $\qbar_\alpha$ constraint (see Equation~\eqref{e:qbarbPoFequiv}). It can be observed that the bPoF-based formulation \eqref{e:bPOFopt1} is equivalent to the $\qbar_\alpha$-based optimization formulation \eqref{e:CVaRopt1} by noting that the bPoF constraint being active implies that the $\qbar_{\alpha}$ constraint is also active, i.e., $\pbuf \left(g(\bd,Z)\right)= 1-\alphaT \implies \qbar_{\alphaT}\left[g(\bd,Z)\right]= t$. 
However, formulation~\eqref{e:bPOFopt1} is useful when considered in the context of interpretability w.r.t.\ the originally intended PoF reliability constraint along with the data-informed conservative buffer provided by bPoF. In engineering applications, the exact failure threshold is often uncertain and chosen by a subject matter expert. Thus, it is beneficial that bPoF can provide a reliability constraint that is robust to uncertain or inexact choices of failure threshold. 

\begin{remark}[Convexity in bPoF-based optimization] \label{r:convexbPOF}
	It can be shown that bPoF\footnote{assuming $\qbar_0\left[g(\bd,Z)\right] < t < \qbar_1\left[g(\bd,Z)\right]$ and $g(\bd,Z)$ is integrable} can be written in the form of a convex optimization problem, similar to $\qbar_\alpha$, as~\cite{norton2019maximization,mafusalov2018buffered}
	\begin{equation} \label{e:bPoF_min}
	\pbuf \left(g(\bd,Z)\right)
	= \min_{\lambda < t}\frac{ \mathbb{E}  \left[ \left[g(\bd,Z) - \lambda \right]^+\right]}{t - \lambda},
	\end{equation}
	where $\bd$ is the given design, $[c]^+=\max\{0,c\}$, and $\lambda$ is an auxiliary variable. Note that the optimal $\lambda^*$ from Equation~\eqref{e:bPoF_min} is the threshold from Equation~\eqref{e:bPoF} that provides $ \mathbb{E}  \left[  g(\bd,Z) \; | \;  g(\bd,Z) \geq \lambda^* \right] = t$.
	
	Using Equation~\eqref{e:bPoF_min}, formulation~\eqref{e:bPOFopt1} can be reduced to an optimization problem involving only expectations as given by
	\begin{equation}
	\label{e:bPoFopt1_gamma}
	\begin{split}
	\min_{\lambda<t,\, \bd\in\mathcal{D}}\quad & \mathbb{E}  \left[f(\bd,Z)\right] \\
	\emph{subject to}\quad & \frac{ \mathbb{E}  \left[ \left[g(\bd,Z) - \lambda \right]^+\right]}{t - \lambda}\le 1-\alphaTrem.
	\end{split}
	\end{equation}
	Note that \eqref{e:bPoFopt1_gamma} can be reformulated by a simple rearrangement of the constraint to become equivalent to the $\qbar_{\alpha}$ constrained problem given by \eqref{e:squantopt1_gammma}. Thus, it is a convex problem when $g(\bd,Z)$ and $f(\bd,Z)$ are convex. The same linearization trick can also be performed as in \eqref{e:squantopt1_gammma_linear}.
	
\end{remark}

\subsubsection{Optimization problem: bPoF as objective}
A bPoF objective provides us with an optimization problem focused on optimal reliability subject to satisfaction of other design metrics. While the use of bPoF as a constraint is equivalent to a $\qbar_\alpha$ constraint, the same can not be said about the case in which bPoF and $\qbar_\alpha$ are used as an objective function.
Consider the PoF minimization problem
\begin{equation} 
\label{e:RBDO2}
\begin{split}
\min_{\bd\in\mathcal{D}}\quad & \pt(g(\bd,Z)) \\
\text{subject to}\quad &  \mathbb{E}  \left[f(\bd,Z)\right]  \le C_{\scriptscriptstyle \text{T}}.
\end{split}
\end{equation}
As mentioned in Remark~\ref{r:discussionRBDO}, PoF is often nonconvex and discontinuous, making gradient calculations ill-posed or unstable. However, \eqref{e:RBDO2} is a desirable formulation if reliability is paramount. The formulation in \eqref{e:RBDO2} defines the situation where given our design performance specifications, characterized by $\mathbb{E}  \left[f(\bd,Z)\right]  \le C_{\scriptscriptstyle \text{T}}$, we desire the most reliable design achievable.

We can consider an alternative to the problem in \eqref{e:RBDO2} using a bPoF objective function as
\begin{equation} 
\label{e:bPOFopt2}
\begin{split}
\min_{\bd\in\mathcal{D}}\quad & \pbuf \left(g(\bd,Z)\right) \\
\text{subject to}\quad &  \mathbb{E}  \left[f(\bd,Z)\right]  \le C_{\scriptscriptstyle \text{T}}.
\end{split}
\end{equation}
Using Equation~\eqref{e:bPoF_min}, the optimization problem in \eqref{e:bPOFopt2} can be rewritten in terms of expectations as
\begin{equation} 
\label{e:bPOFopt3}
\begin{split}
\min_{\lambda < t, \ \bd\in\mathcal{D}} & \frac{ \mathbb{E}  \left[ \left[g(\bd,Z) - \lambda \right]^+\right]}{t - \lambda},\\
\text{subject to}\quad &  \mathbb{E}  \left[f(\bd,Z)\right]  \le C_{\scriptscriptstyle \text{T}}.
\end{split}
\end{equation}
This allows one to minimize bPoF, a conservative upper bound of PoF, subject to constraints on costs or other performance measures, while maintaining convexity and other mathematical properties of the limit state function.

\subsubsection{Discussion on bPoF-based optimization}
There are several advantages of using the bPoF-based optimization problem described by Equation~\eqref{e:bPOFopt1} as compared to the RBDO problem. First, the bPoF-based optimization problem leads to a data-informed conservative design as compared to RBDO. This data-informed conservativeness of the bPoF-based optimal design is more desirable and resilient because it takes into account the magnitude of failure (or the tail of the distribution) that guards against more serious catastrophic failures as highlighted later in Remark~\ref{r:conserve}. Second, the bPoF-based optimization problem preserves convexity if the underlying limit state function is convex as noted in Remark~\ref{r:convexbPOF} (as compared to PoF, which does not preserve convexity). This leads to well-behaved convex optimization problems even for the risk-based formulation and pushes us to pay more attention to devising limit state functions that are convex or nearly convex. Additionally, if $\pi$ is independent of $\bd$, the same linear reformulation trick used with $\qbar_\alpha$ constraints (see \eqref{e:squantopt1_gammma_linear}) can be used to transform the objective in \eqref{e:bPOFopt3} into a linear function with additional linear constraints and auxiliary variables, offering similar advantages as noted in Section~\ref*{s:3}.\ref*{s:SQopt}.\ref{s:cvarDisc}. Third, under certain conditions, it is possible to calculate (quasi)-gradients for bPoF~\cite{zhang2019derivatives}. Note that non-smoothness can also be avoided by using smoothed versions given by \cite{kouri2019higher,basova2011computational}. 

\section{Structural Design: Short Column Problem With Convex Reformulation}\label{s:4}
In this section, we use the short column structural design problem that has been widely used in the RBDO community~\cite{aoues2010benchmark,bichon2013efficient} to compare some of the properties of PoF- and bPoF-based optimization formulations. We show a convex reformulation of the objective and limit state functions leading to strongly certifiable risk-based designs for bPoF-based optimization. 

\subsection{Short column problem description}
The problem consists of designing a short column with rectangular cross-section of dimensions $w$ and $h$, subjected to uncertain loads (axial force $F$ and bending moment $M$). The yield stress of the material, $Y$, is also considered to be uncertain. The random variables are $Z = [ F, \ M, \ Y]^\top $ with a joint distribution $\pi$. Table~\ref{t:randvarSC} describes the random variables used in the short column design. The correlation coefficient between $F$ and $M$ is 0.5. The design variables, $\bd = [ w, \ h]^\top$, are the width and depth of the cross-section as shown in Table~\ref{t:desvarSC}. The objective function is the cross-sectional area given by $wh$. Along with a failure threshold $t=1$, the limit state function is defined as
\begin{equation}
g(\bd,\bz)=\frac{4M}{wh^2Y}+\frac{F^2}{w^2h^2Y^2}.
\end{equation}
\begin{table}[!htb]
	\centering
	\caption{Random variables used in the short column application.}
	\label{t:randvarSC}
	\begin{tabular}{C{2.5cm}L{2cm}C{2cm}C{2cm}C{3cm}}
		\hline
		Random variable & Units & Distribution & Mean & Standard deviation \\ 
		\hline
		$F$ & kN & Normal & 500 & 100  \\
		$M$ & kNm & Normal & 2000 & 400 \\
		$Y$ & MPa & Log-normal & 5 & 0.5  \\
		\hline
	\end{tabular}
\end{table}
\begin{table}[!htb]
	\centering
	\caption{Design variables used in the short column application.}
	\label{t:desvarSC}
	\begin{tabular}{ccc}
		\hline
		Design variable & Lower bound (m) & Upper bound (m) \\ 
		\hline
		$w$ & 5 & 15    \\
		$h$ & 15 & 25 \\
		\hline
	\end{tabular}
\end{table}

\subsection{Optimization problem formulations}
This section provides the optimization formulations based on PoF and bPoF for the short column structural design. We show a convex reformulation of the bPoF-based optimization short column problem to emphasize the specific advantage of bPoF risk measure making it strongly certifiable. For each case, we solve multiple optimization problems each with a different fixed value of $1-\alphaT$, where $\alphaT$ is the desired reliability.

\subsubsection{Short column RBDO}
The RBDO problem is given by
\begin{equation}
\label{e:RBDO_SC}
\begin{split}
\min_{w,h}\quad & wh \\
\text{subject to}\quad & \pt(g(\bd,Z))\le 1-\alphaT,\\
& \ell_w \leq w \leq u_w, \\
& \ell_h \leq h \leq u_h,
\end{split}
\end{equation}
where $(\ell_w, \ell_h,u_w,u_h)$ denote the lower and upper bounds on $w$ and $h$ as defined in Table~\ref{t:desvarSC}. 

\subsubsection{Short column bPoF-constrained CRiBDO with convex reformulation}
The bPoF-based optimization problem for the short column design is
\begin{equation}
\label{e:short_col_cvar0}
\begin{split}
\min_{\lambda<t, w,h}\quad & wh \\
\text{subject to}\quad & \frac{ \mathbb{E} \left[ \left[ \frac{4M}{wh^2Y}+\frac{F^2}{w^2h^2Y^2} - \lambda \right]^+\right]}{t-\lambda}  \le 1-\alphaT, \\
& \ell_w \leq w \leq u_w, \\
& \ell_h \leq h \leq u_h,
\end{split}
\end{equation}
where we used \eqref{e:bPoFopt1_gamma}. We first show that the optimization problem with a bPoF constraint can be formulated as a convex optimization problem in this case. This will permit us to take advantage of convex optimization solvers and offer guarantees for the optimization problem. To reformulate this as a more manageable convex optimization problem\footnote{We were required to find a formulation that was recognized as convex by modeling language CVXpy with convex solver MOSEK}, we rearrange the first constraint in \eqref{e:short_col_cvar0} to achieve an equivalent form\footnote{Note that this is in the same form as a $\qbar_{\alphaT}$ constraint from (\ref{e:squantopt1_gammma}).},
\begin{equation}
\label{e:short_col_cvar1}
\lambda+ \frac{1}{(1-\alphaT)} \mathbb{E} \left[ \left[ \frac{4M}{wh^2Y}+\frac{F^2}{w^2h^2Y^2} - \lambda \right]^+\right]  \le t.
\end{equation}

Next, we note that both $w\text{ and }h$ are nonnegative and thus we can make the change of variable $w=e^{x_1} , h = e^{x_2}$ with $x_1,x_2 \in \mathbb{R}$. Then the limit state function becomes $g(x_1,x_2,\bz) = \frac{4M}{Y}e^{-x_1-2x_2}+\frac{F^2}{Y^2}e^{-2x_1-2x_2}$ and objective function becomes $wh = e^{x_1 + x_2}$. These are both convex functions in the new design variables $(x_1,x_2)$. Thus, using the change of variable and Equation~\eqref{e:short_col_cvar1} we can reformulate \eqref{e:short_col_cvar0} as
\begin{equation}
\label{e:short_col_cvar2}
\begin{split}
\min_{\lambda<t, x_1,x_2}\quad & e^{x_1 + x_2} \\
\text{subject to}\quad &  \lambda+ \frac{1}{(1-\alphaT)} \mathbb{E} \left[\left[ \frac{4M}{Y}e^{-x_1-2x_2}+\frac{F^2}{Y^2}e^{-2x_1-2x_2} - \lambda \right]^+\right]  \le t, \\
&\ln \ell_w \leq x_1 \leq \ln u_w, \\
& \ln \ell_h \leq x_2 \leq \ln u_h.
\end{split}
\end{equation}
Furthermore, since the distribution of the random variables are independent of the design variables, we can empirically estimate the expectation in the constraint for any design by using a fixed set of $m$ samples $\{\bz_1,\dots,\bz_m\}$ that are sampled \textit{a priori} from the distribution $\pi$. This gives us a convex sample-average-approximation optimization problem for the given set of $m$ samples as
\begin{equation}
\label{e:short_col_cvar3}
\begin{split}
\min_{\lambda<t, x_1,x_2}\quad & e^{x_1 + x_2} \\
\text{subject to}\quad &  \lambda+ \frac{1}{m(1-\alphaT)} \sum_{i=1}^m \left[ \frac{4M_i}{Y_i}e^{-x_1-2x_2}+\frac{F_i^2}{Y_i^2}e^{-2x_1-2x_2} - \lambda \right]^+  \le t, \\
&\ln \ell_w \leq x_1 \leq \ln u_w \\
& \ln \ell_h \leq x_2 \leq \ln u_h .
\end{split}
\end{equation}
Note that the constraint is convex since $[\cdot]^+$ is a convex function and preserves the convexity of the limit state function. Thus, the formulation in \eqref{e:short_col_cvar3} leads to a strongly certifiable formulation since it satisfies both conditions of certifiability.

\subsection{Experimental comparison between RBDO and bPoF-based CRiBDO}
We now compare the behavior of the RBDO formulation and the bPoF-constrained CRiBDO formulation. We solve the RBDO problem with various values of $\alphaT$ using the gradient-free COBYLA optimizer~\cite{powell1994direct}.  We estimate PoF in the RBDO problem by iteratively adding samples until the MC error reaches below 1\%. We solve the bPoF-constrained problem (Equation~\eqref{e:short_col_cvar3}) with various values of $\alphaT$, utilizing the convex optimization solver MOSEK with CVXPY~\cite{diamond2016cvxpy} as our modeling interface to the solver. 

\begin{figure}[!htb]
	\centering
	\subfigure[]{\includegraphics[width=8.2cm,page=2]{Results_all_SC.pdf}}\hfill
	\subfigure[]{\includegraphics[width=8.2cm,page=3]{Results_all_SC.pdf}}\hfill
	\subfigure[]{\includegraphics[width=8.2cm,page=1]{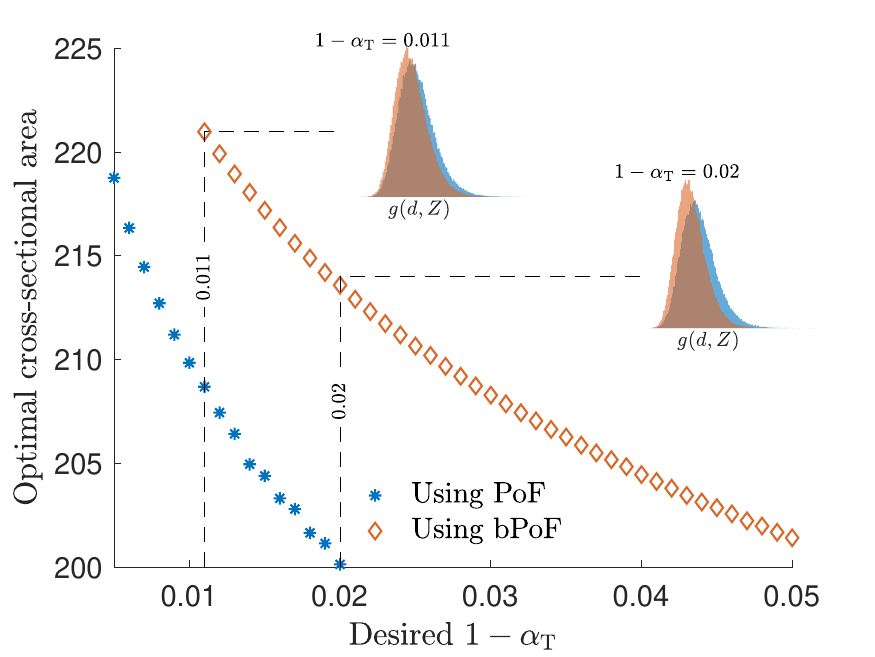}}\hfill
	\subfigure[]{\includegraphics[width=8.2cm,page=9]{Results_all_SC.pdf}}
	\caption{Optimal cross-sectional area and designs obtained for different levels of target $1-\alphaT$.}
	\label{fig:SCarea}
\end{figure}

We illustrate that, as hypothesized, PoF and bPoF are indeed natural counterparts with both measuring notions of failure risk. Each point in Figure~\ref{fig:SCarea}(a) corresponds to an optimal design achieved by solving RBDO (``Using PoF'') or bPoF-constrained optimization for some $\alphaT$ (``Using bPoF''). The x-axis provides the estimated value of PoF of the design, estimated using a separate evaluation sample of size $5 \times 10^6$. The y-axis provides the cross-section area of the optimal design. Figure~\ref{fig:SCarea}(b) is similar, except with the value of bPoF (estimated using $5 \times 10^6$ samples) as the x-axis. We see that RBDO and bPoF-based CRiBDO achieve nearly identical frontiers of cross-section area for the optimal designs with similar PoF or bPoF.  The optimal designs similar to the ones obtained using RBDO can be obtained by setting higher values for $1-\alphaT$ in the bPoF constraint as compared to PoF constraint, and vice versa. This highlights the fact that bPoF is indeed a natural replacement for PoF, achieving similar goals and controlling failure probability. All the optimal designs for the short column problem have $h\approx25$ as shown in Figure~\ref{fig:SCarea}(d).
However, an advantage of the bPoF-based CRiBDO formulation is that while both achieve similar design frontiers, we have a guarantee that the design given by the bPoF formulation is \textit{globally optimal} for the given sample-average-approximation problem in Equation~\eqref{e:short_col_cvar3} since it is a convex problem. We have no such guarantee for the designs given by RBDO. Thus, even if similar designs are achieved by RBDO and bPoF-constrained optimization, we have additional guarantees about the quality of the bPoF-based design due to the underlying convexity. 

\begin{figure}[!htb]
	\centering
	\subfigure[RBDO designs]{\includegraphics[width=0.5\textwidth,page=6]{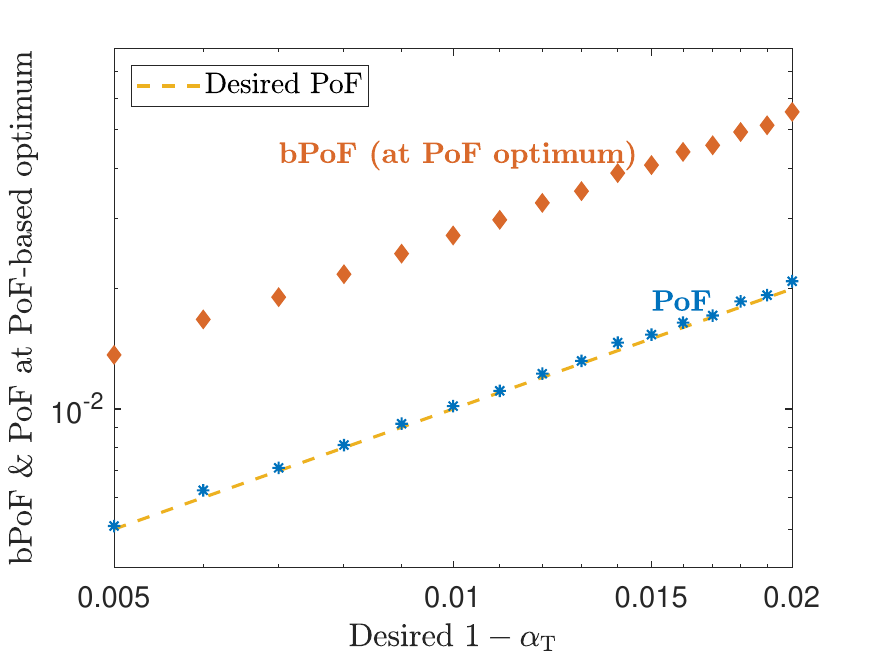}}\hfill
	\subfigure[bPoF-based CRiBDO designs]{\includegraphics[width=0.5\textwidth,page=7]{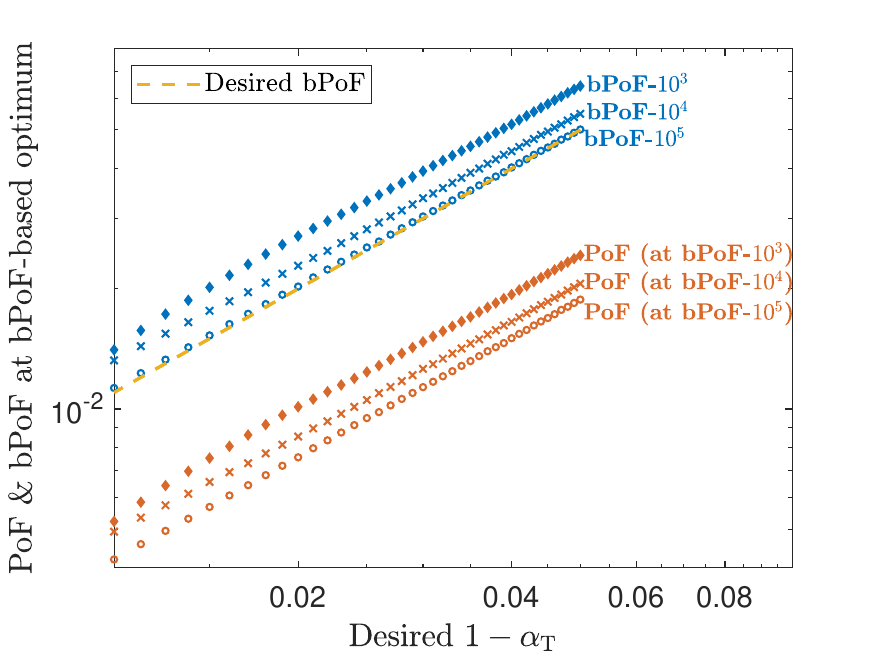}}
	\caption{Comparing (a) bPoF estimates at the optimal designs obtained through RBDO, and (b) PoF estimates at the optimal designs obtained through bPoF-based CRiBDO using different \textit{a priori} sample sizes\protect\footnotemark\ for different desired $1-\alphaT$ in log-log scale.}
	\label{fig:SCPoFbPoF}
\end{figure}

Another advantage is the conservative nature of bPoF as compared to PoF. When formulated with identical levels of $\alphaT$ in the constraint, the bPoF-based optimization achieves a more conservative design than RBDO. The conservative nature of bPoF for the same desired reliability $\alphaT$ can be seen from Figure~\ref{fig:SCPoFbPoF} and is directly reflected in more conservative optimal designs with larger cross-sectional areas as seen in Figure~\ref{fig:SCarea}(c). This type of probabilistically derived conservativeness can be seen as desirable as highlighted in Remark~\ref{r:conserve} below.

Figure~\ref{fig:SCPoFbPoF} compares different levels of desired $1-\alphaT$ versus the estimated PoF or bPoF for the optimal designs obtained through PoF- and bPoF-based optimization. For these plots, we use $5\times10^6$ samples\footnote{These samples are not used in the optimization, but only to estimate at the optimal design after the optimization is completed.} to get accurate estimates of the PoF or bPoF at the optimum. Figure~\ref{fig:SCPoFbPoF}(a) shows the desired $1-\alphaT$ and the PoF/bPoF for the optimum design obtained using the RBDO problem. We can see that since the MC error for PoF estimate in the RBDO problem was always ensured to be below 1\%, the desired PoF and the PoF at the optimum overlap. The figure also emphasizes the conservative property of bPoF for the same desired $1-\alphaT$. 

A key observation is illustrated by Figure~\ref{fig:SCPoFbPoF}(b), which compares the results for optimal designs obtained using bPoF-based CRiBDO for different \textit{a priori} sample sizes, i.e. the value of $m$ in Equation~\eqref{e:short_col_cvar3}. We make this comparison to analyze the effect of fixing the sample set for all optimization iterations before starting the optimization, which is required to obtain the convex optimization formulation shown in Equation~\eqref{e:short_col_cvar3}. We can see that for lower sample sizes of $10^3$ and $10^4$, the bPoF at the optimum and the desired bPoF do not overlap reflecting inaccurate MC estimates of bPoF. However, it should be noted that the bPoF formulation is still effective in controlling the PoF, even when sample size is small. In other words, even when small number of samples are used within the optimization, the nature of bPoF yields an optimal design with desirable conservativeness and thus, an acceptably low PoF. Additionally, even when formulated with a small sample size, the bPoF-based convex optimization problem is still considerably stable leading to good optimal designs. One of the primary drawbacks of RBDO is the potential fragility of the optimization, particularly when sample sizes are small, where the estimates of PoF and/or gradients (if a gradient-based solver is used) are unstable and produce poor or inconsistent optimization results. The bPoF formulation does not seem to suffer in the same way for the short column design as illustrated here. 
\footnotetext{The \textit{a priori} sample sizes $m$ used for the bPoF CRiBDO are indicated in the legend using bPoF-$m$.}

\begin{remark}[Desirable data-informed conservativeness] \label{r:conserve}
	We take a closer look at the conservativeness induced by the bPoF-based CRiBDO for the same desired reliability level $\alphaT$ as compared to PoF-based optimization (as shown in Figures~\ref{fig:SCarea}(c) and \ref{fig:SCPoFbPoF}) and why this type of conservativeness would be desirable and resilient. There are other ways of introducing conservativeness, such as safety factors, basis values, and stricter reliability levels. Typically, using safety factors leads to overly conservative designs. When inappropriate safety factor values are used, the deterministic optimization setup could also potentially lead to unreliable designs. This is because converting to deterministic optimization using just safety factors (or basis values) to account for the uncertainty in the system does not take into account the distribution of the limit state function and lacks sufficient information to make good decisions. These well-known issues with safety-factor- and basis-values-based deterministic optimization formulations have progressively led us to consider risk-based optimization under uncertainty.
	Another way to introduce conservativeness in risk-based optimization is by using lower values of $1-\alphaT$, which leads to stricter reliability constraints. However, this will just lead to overly reliable designs without any information about the distribution of the limit state function. 
	
	The bPof-based CRiBDO can be seen as a better way to induce conservativeness because it encodes more information about the underlying limit state function through the data on the magnitude of failure (as seen from Equations~\eqref{e:bPoF} and \eqref{e:bPOFconserve}). In Section~\ref{s:5}, we highlight a similar observation on conservativeness for $\qbar_{\alpha}$-based CRiBDO through the thermal design problem. These CRiBDO formulations lead to a probabilistic data-informed way of achieving a conservative design, which can be seen as more desirable in practice.
\end{remark}

\section{Thermal Design: Cooling Fin Problem}\label{s:5}
In this section, we compare the properties of PoF- and $\qbar_\alpha$-based optimization formulations for the thermal design of a cooling fin problem. A cooling fin is used in many engineering applications to dissipate heat. This application also highlights the usefulness of CRiBDO even without the second certifiability condition of convexity. This problem considers a discretized partial differential equation, where closed-form solutions (such as in the previous short column) are not available, and only computer simulations are available.

\subsection{Cooling fin model description}
\begin{figure}[!h]
	\centering
	\includegraphics[width=0.5\textwidth,trim = {0.5cm 0.5cm 0.5cm 0.5cm},clip]{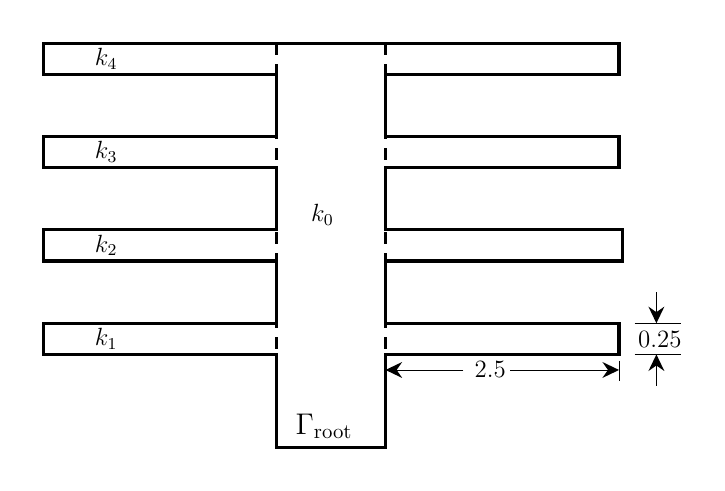}
	\caption{Fin geometry and model parameters.}
	\label{fig:finGeom}
\end{figure}

We consider a cooling fin with fixed geometry as shown in Figure~\ref{fig:finGeom}, consisting of a vertical post with horizontal fins attached.  We briefly review the problem here and refer to \cite{CPrudhomme_DVRovas_KVeroy_LMachiels_YMaday_ATPatera_GTurinici_2002a} for more details. The fin array consists of four horizontal sub-fins with width 2.5 and thickness 0.25, as well as a fin post with unit width and height of four. The thermal design is parametrized by the fin conductivities $k_i, i =1, \ldots,4$ and the post conductivity $k_0$, as well as the Biot number $Bi$, which is a non-dimensionalized heat transfer coefficient for thermal transfer from the fins to the surrounding air. 
The design variables, $\bd = [k_1, k_2, k_3, k_4]$ are the thermal conductivities of the four fins as shown in Table~\ref{t:desvarTF}. The post conductivity is $k_0=5$ and the Biot number is $Bi=0.5$. We introduce manufacturing and operational uncertainties in all the parameters through the random variable $Z = [ \xi_0, \ \xi_1, \ \xi_2, \ \xi_3, \ \xi_4,  \xi_{Bi}]^\top$ with a joint distribution $\pi$ given in Table~\ref{t:randvarTF}. The random variables $\xi_i$ model the additive uncertainty for the respective thermal conductivities $k_i, i=0,\dots,4$ and $\xi_{Bi}$ models the additive uncertainty for the Biot number $Bi$.
The system is governed by Poisson's equation in two spatial dimensions denoted by $\boldsymbol{x}$ whose solution is the temperature field $y(\boldsymbol{x}, \bd, Z)$. The PDE is semi-discretized with the finite element method and yields a system with $4,760$ degrees of freedom. 

The fin conducts heat away from the root $\Gamma_{\text{root}}$, so the lower the root temperature, the more effective the cooling fin. Thus, our objective function depends on the measure of the average temperature at the root, i.e., 
\begin{equation}
\mathcal{Y}(\bd,Z) = \int_{\Gamma_{\text{root}}} y(\boldsymbol{x}, \bd, Z) \text{d} x .
\end{equation}
We also include a quantity proportional to the cost of the material based on the area and material thermal conductivity in the objective function as shown in Section~\ref*{s:5}.\ref{s:TFoptform}. The limit state function for the cooling fin problem is based on the maximum temperature and is defined as
\begin{equation}
g(\bd,Z)=\max_{x} y(\boldsymbol{x}, \bd, Z).
\end{equation}
We choose $t=0.35$ as the constraint on the limit state function to define the maximum allowable temperature of the system.

\begin{table}[!htb]
	\centering
	\caption{Design variables used in the cooling fin application.}
	\label{t:desvarTF}
	\begin{tabular}{ccc}
		\hline
		Design variable & Lower bound & Upper bound \\ 
		\hline
		$k_i, \ i=1, \ldots 4$ & 1 & 10 \\
		\hline
	\end{tabular}
\end{table}

\begin{table}[!htb]
	\centering
	\caption{Random variables used in the cooling fin application.}
	\label{t:randvarTF}
	\begin{tabular}{C{2.5cm}C{6cm}C{1.5cm}C{3cm}}
		\hline
		Random variable & Distribution & Mean $\mu$ & Standard deviation $\sigma$ \\ 
		\hline
		$\xi_i, \ i=0, \ldots 4$ & \multirow{2}{6cm}{truncated normal ($[\mu-4\sigma,\mu+2\sigma]$)} & 0 & 0.1  \\
		$\xi_{Bi}$ &  & 0 & 0.02 \\
		\hline
	\end{tabular}
\end{table}

\subsection{Optimization problem formulations}\label{s:TFoptform}
This section provides the optimization formulations based on PoF and $\qbar_\alpha$ for the cooling fin thermal design.

\subsubsection{Cooling fin RBDO}
The RBDO problem is given by
\begin{equation}
\label{e:RBDO_TF1}
\begin{split}
\min_{\bd\in\mathcal{D}}\quad & \mathcal{Y}(\bd,\mu_Z) + \frac{\sum_{i=0}^{4}A_ik_i}{5A_0 + \sum_{i=1}^{4}10A_i}\\
\text{subject to}\quad & \pt(g(\bd,Z))\le 1-\alphaT, 
\end{split}
\end{equation}
where $A_i$ denotes the area for the material with thermal conductivity of $k_i, i=0,\dots,4$ and $A_ik_i$ represents a quantity proportional to the cost of the material. Here, the fin post area is $A_0=4$ and sub-fin areas are $A_i = 1.25, i=1\dots,4$. The cost part is normalized by the maximum proportionate cost. We use two different values of $1-\alphaT \in \{0.001, 0.05\}$.

\subsubsection{Cooling fin $\qbar_\alpha$-constrained CRiBDO}
The $\qbar_\alpha$-constrained CRiBDO formulation for the cooling fin design is
\begin{equation}
\label{e:TF_cvar0}
\begin{split}
\min_{\bd\in\mathcal{D}}\quad & \mathcal{Y}(\bd,\mu_Z) + \frac{\sum_{i=0}^{4}A_ik_i}{5A_0 +\sum_{i=1}^{4}10A_i}\\
\text{subject to}\quad & \qbar_{\alphaT} \left[g(\bd,Z)\right]\le t,
\end{split}
\end{equation}
where we find the optimal designs for two different values of $1-\alphaT \in \{0.001, 0.05\}$. In this case, the underlying limit state function is not known to be convex making the CRiBDO formulation \eqref{e:TF_cvar0} certifiable in one condition, which is the data-informed conservativeness.

\subsection{Experimental comparison between RBDO and $\qbar_\alpha$-constrained CRiBDO}
We compare the optimal results obtained through RBDO and $\qbar_\alpha$-constrained CRiBDO formulations under the same $\alphaT$ values. We solve the RBDO and CRiBDO problems using the gradient-free COBYLA optimizer. We estimate the PoF in each RBDO iteration by iteratively adding samples until the MC error reaches below 1\% with the maximum number of samples capped at $10^5$. We estimate the $\qbar_\alpha$ in each CRiBDO iteration by using $10^5$ MC samples.

\begin{figure}[!htb]
	\centering
	\hspace*{\fill}
	\subfigure[$1-\alphaT=0.001$]{\includegraphics[width=0.45\textwidth, trim = {0.5cm 0cm 0.8cm 0.5cm},clip,page=1]{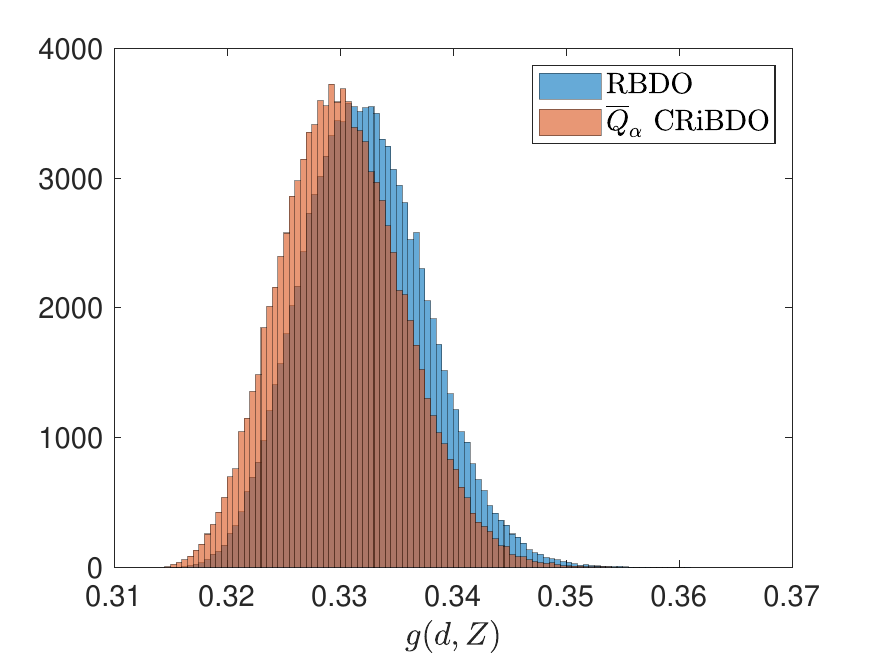}}\hfill
	\subfigure[$1-\alphaT=0.05$]{\includegraphics[width=0.45\textwidth, trim = {0.5cm 0cm 0.8cm 0.5cm},clip,page=2]{Plot_RBDOCVaRoptDesigns_All.pdf}}	
	\hspace*{\fill}
	\caption{Histograms comparing limit state function distributions for optimal designs obtained through different optimization formulations.}
	\label{fig:TF_dOptDist}
\end{figure}

\begin{figure}[!bht]
	\centering
	\subfigure[RBDO ($1-\alphaT=0.001$)]{\includegraphics[width=0.4\textwidth, trim = {0.5cm 0cm 0.8cm 0.5cm},clip,page=3]{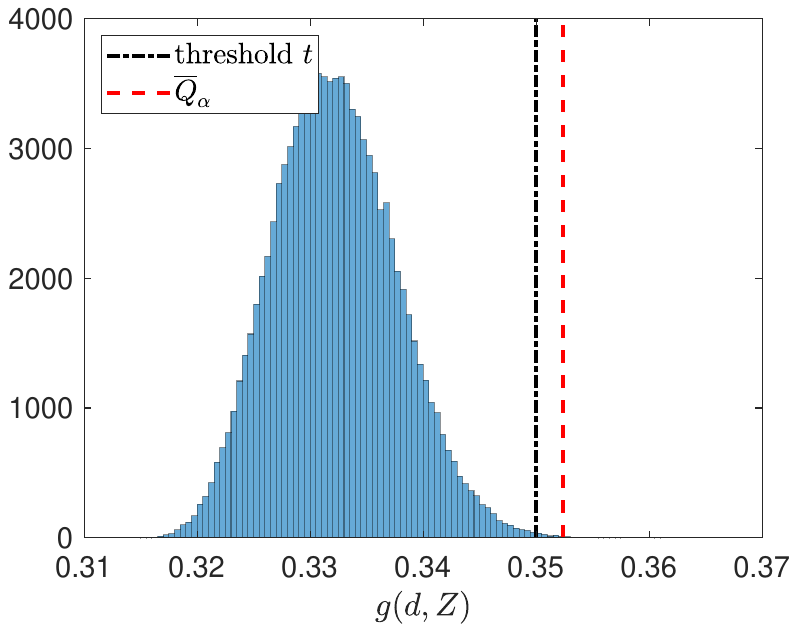}}\hspace{2em}
	\subfigure[RBDO ($1-\alphaT=0.05$)]{\includegraphics[width=0.4\textwidth, trim = {0.5cm 0cm 0.8cm 0.5cm},clip,page=4]{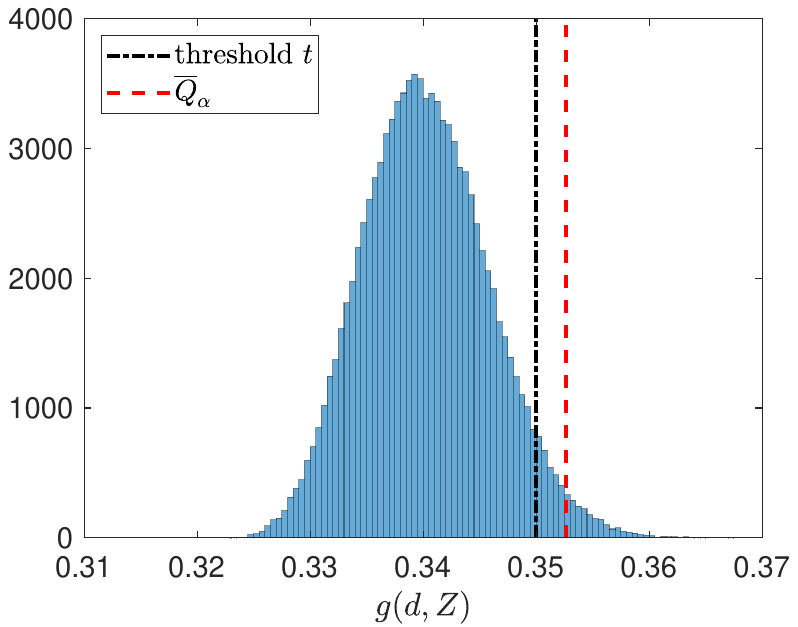}}
	\subfigure[$\qbar_\alpha$ CRiBDO ($1-\alphaT=0.001$)]{\includegraphics[width=0.4\textwidth, trim = {0.5cm 0cm 0.8cm 0.5cm},clip,page=5]{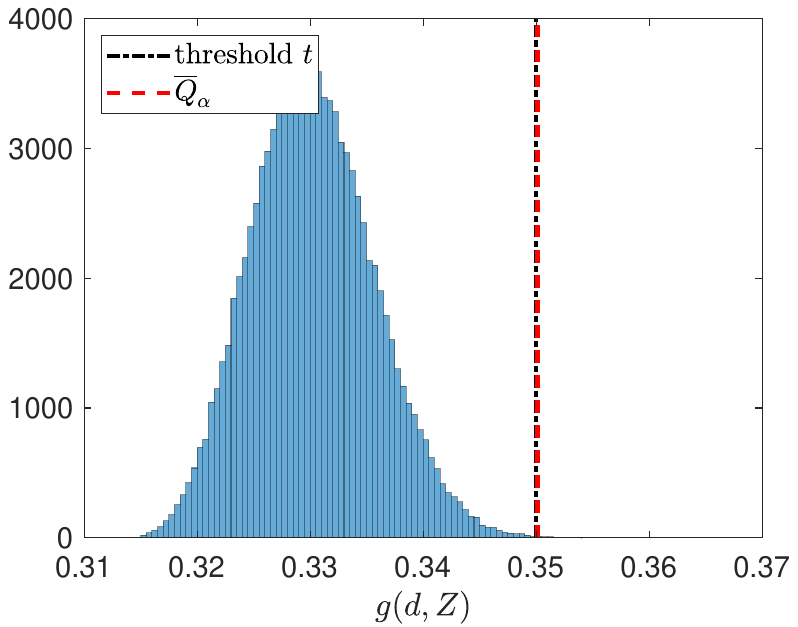}}\hspace{2em}
	\subfigure[$\qbar_\alpha$ CRiBDO ($1-\alphaT=0.05$)]{\includegraphics[width=0.4\textwidth, trim = {0.5cm 0cm 0.8cm 0.5cm},clip,page=6]{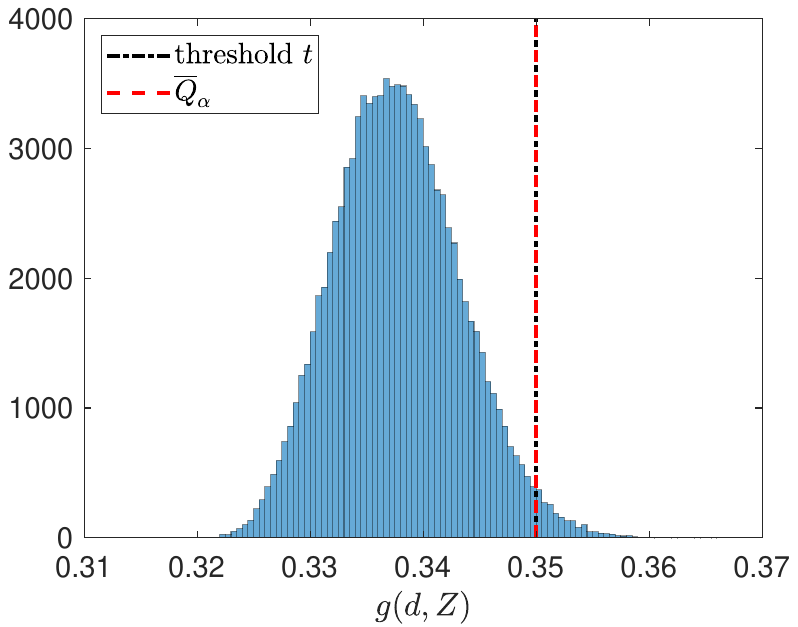}}
	\caption{Comparing specified thresholds and $\qbar_\alpha$ for optimal designs obtained through different optimization formulations.}
	\label{fig:TF_dOptThresh}
\end{figure}

Table~\ref{t:resultsTF} shows the optimal designs obtained from the different optimization formulations. We start the optimization with an initial design that is feasible for all the optimization formulations. The optimal designs obtained using $\qbar_\alpha$-constrained CRiBDO for a given $\alphaT$ are more conservative than the RBDO designs. This highlights one of the major advantages of using $\qbar_\alpha$-constrained CRiBDO that certifies designs through the data-informed conservativeness. The conservative nature of the $\qbar_\alpha$-constrained CRiBDO can be clearly seen by comparing the limit state function distributions at the optimal designs as shown in Figure~\ref{fig:TF_dOptDist}. As discussed in Remark~\ref{r:conserve}, this conservativeness is desirable and required to prevent catastrophic failures. Figure~\ref{fig:TF_dOptThresh} compares the specified hard thresholds with the $\qbar_\alpha$ for the different optimal designs to further highlight the fact that $\qbar_\alpha$ considering the magnitude of failure and not using hard thresholding leads to appropriately conservative designs. The data-informed nature of the conservativeness is a significant advantage since the magnitude of conservativeness induced automatically changes according to the data from the underlying limit state function distribution for a particular design, i.e., $\qbar_\alpha$ is more conservative only when it is required as dictated by the underlying distribution. We explicitly show the data-informed nature of conservativeness in the next section.
\begin{table}[!htb]
	\centering
	\caption{Optimal designs obtained from RBDO and $\qbar_\alpha$-constrained CRiBDO.}
	\label{t:resultsTF}
	\begin{tabular}{C{2.6cm}C{3.3cm}C{2.1cm}C{2cm}C{2.1cm}C{2cm}}
		\hline
		\multirow{2}{2.6cm}{Design variable/ Output statistic} & \multirow{2}{1cm}{Initial design} & \multicolumn{2}{c}{RBDO} & \multicolumn{2}{c}{$\qbar_\alpha$-constrained CRiBDO} \\
		\cline{3-4} \cline{5-6}
		& & $1-\alphaT=0.001$ & $1-\alphaT=0.05$ & $1-\alphaT=0.001$ & $1-\alphaT=0.05$ \\
		\hline
		$k_1$ & 5.5 & 4.154 & 3.3271 & 4.5175  & 3.6468 \\
		$k_2$ & 5.5 & 4.1132 & 3.2015 & 3.965 & 3.1974 \\
		$k_3$ & 5.5 & 1 & 1.0053 & 1 & 1 \\
		$k_4$ & 5.5 & 1 & 1 & 1 & 1 \\
		Objective function & 1.0141 & 0.7649 & 0.7094 & 0.781 & 0.725 \\
		PoF & 0 & $0.001$ &  0.0486 & $\mathit{4.3\times 10^{-4}}$ & \textit{0.0195} \\
		$\qbar_{\alphaT}$ & \parbox[t]{3.3cm}{0.3376 ($1-\alphaT=0.001$)\\0.3301 ($1-\alphaT=0.05$)} & \textit{0.3524} & \textit{0.3526} & 0.35 & 0.35 \\
		\hline
	\end{tabular}
\end{table}

\subsection{Effect of input uncertainties on data-informed conservativeness}\label{s:TF_conserve}
\begin{figure}[!htb]
	\centering
	\includegraphics[width=8cm, trim = {0.5cm 0cm 0.8cm 0.5cm},clip,page=1]{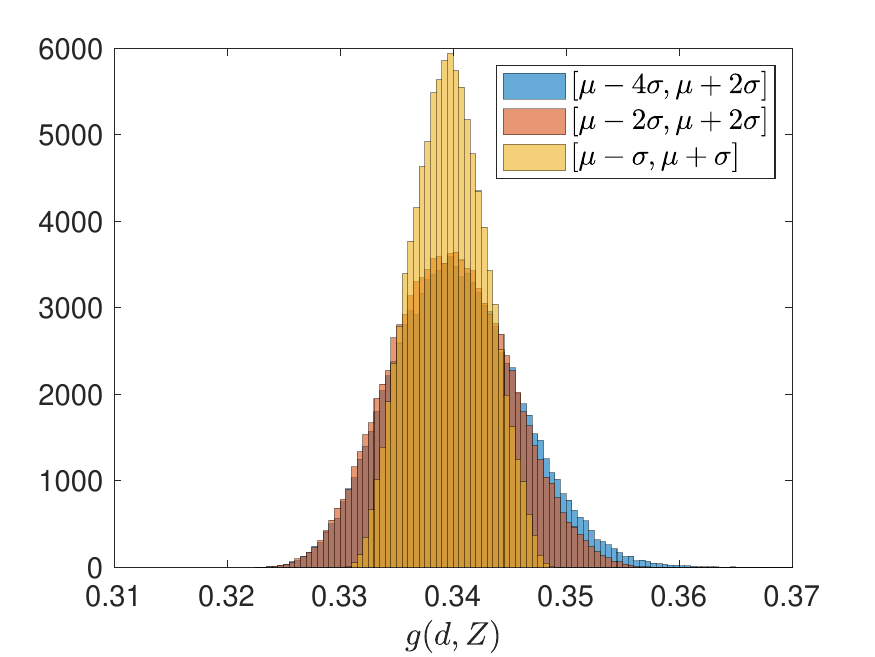}	
	\caption{Histograms for different limit state function distributions generated through modifying input uncertainty truncation range for a fixed design $\bd = [3.3271,3.2015,1.0053,1]$.}
	\label{fig:TF_DiffLimDist}
\end{figure}

In this section, we demonstrate the data-informed nature of the conservativeness induced by $\qbar_{\alpha}$ as described in Remark~\ref{r:conserve}.  The magnitude of conservativeness is naturally adjusted for different limit state function distributions. Since, $\qbar_{\alpha}$ and $Q_{\alpha}$ are natural counterparts, we quantify the magnitude of conservativeness by the percentage difference $(\qbar_{\alpha}-Q_\alpha)/Q_{\alpha}\%$. We fix the design and the $\alpha$ value for the comparison in this section. We use the optimal design obtained by RBDO with $1-\alphaT=0.05$ given in Table~\ref{t:resultsTF} as the fixed design $\bd = [3.3271,3.2015,1.0053,1]$. Different limit state function distributions are generated by changing the input uncertainties through modifying the truncation range for the distribution given in Table~\ref{t:randvarTF}. Figure~\ref{fig:TF_DiffLimDist} shows histograms of the different limit state function distributions for the fixed design. Three different truncation ranges for the input random variables are indicated in the Figure~\ref{fig:TF_DiffLimDist} legend. The distribution of the limit state function is obtained by first sampling the input random variables $Z$ from their given distribution and then running the simulation to obtain realizations of $g(\bd,Z)$. All the results are generated using $10^5$ MC samples. 

\begin{figure}[!htb]
	\centering
	\subfigure[{$[\mu-4\sigma,\mu+2\sigma]$}]{\includegraphics[width=5.5cm, trim = {0.5cm 0cm 0.8cm 0.5cm},clip,page=2]{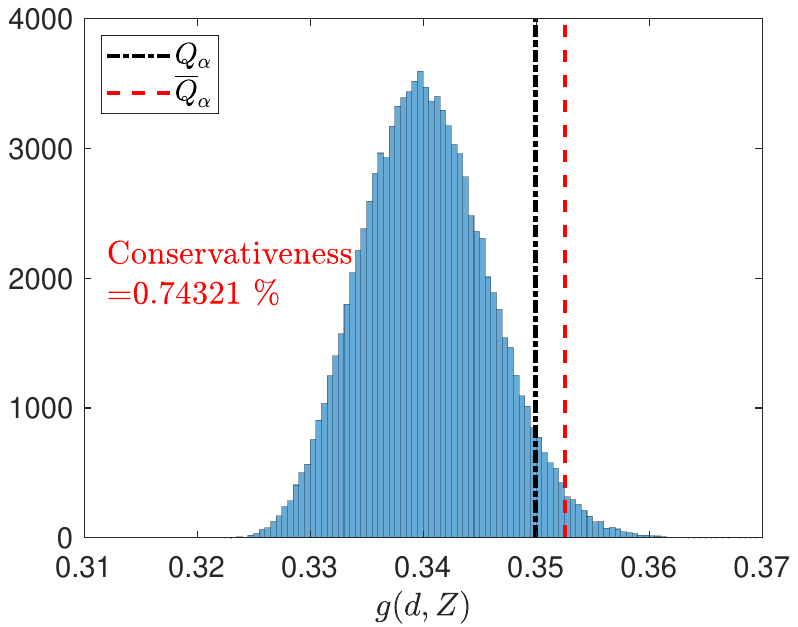}}\hfill
	\subfigure[{$[\mu-2\sigma,\mu+2\sigma]$}]{\includegraphics[width=5.5cm, trim = {0.5cm 0cm 0.8cm 0.5cm},clip,page=3]{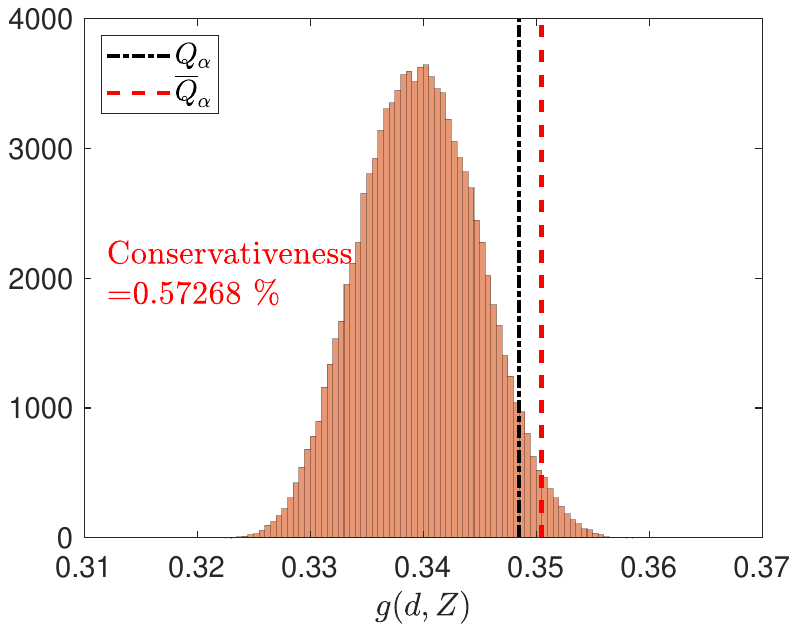}}\hfill
	\subfigure[{$[\mu-\sigma,\mu+\sigma]$}]{\includegraphics[width=5.5cm, trim = {0cm 0cm 0.8cm 0.5cm},clip,page=4]{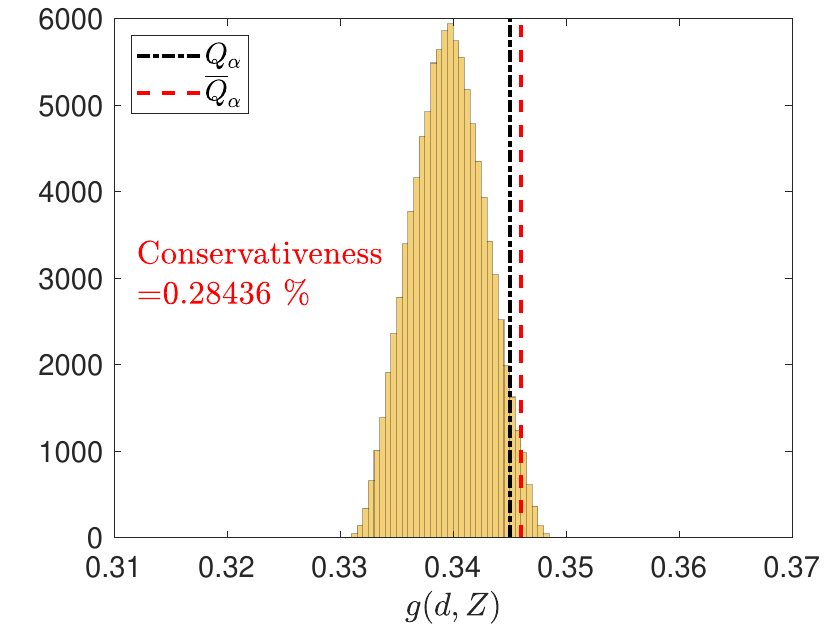}}
	\caption{Conservativeness induced by $\qbar_{\alpha}$ compared to $Q_{\alpha}$ for different limit state function distributions generated through modifying input uncertainty truncation range for a fixed design $\bd = [3.3271,3.2015,1.0053,1]$ and $1-\alpha=0.05$.}
	\label{fig:TF_cvarConserve}
\end{figure}

Figure~\ref{fig:TF_cvarConserve} shows the different levels of conservativeness of $\qbar_{\alpha}$ when compared to $Q_{\alpha}$ for different limit state function distributions. We can see that $\qbar_{\alpha}$ is always conservative when compared to $Q_{\alpha}$. Furthermore, it can be seen that the magnitude of conservativeness depends on the distribution, which exemplifies the data-informed nature of the induced conservativeness, i.e., $\qbar_{\alpha}$ is as conservative as required by the underlying distribution. Specifically, the magnitude of conservativeness for superquantile is higher for the fatter-tailed distribution as seen in Figure~\ref{fig:TF_cvarConserve} (a). For fat-tailed distributions, i.e., distributions with significant tail risk, superquantile provides additional conservativeness by nature of being a tail-integral.

\section{Concluding Remarks}\label{s:6}
In this work, we propose two certifiability conditions that lead to certifiable risk-based design optimization (CRiBDO): (a) data-informed conservativeness: the resulting designs should be certifiably risk-averse against near-failure and catastrophic failure events, and (b) optimization convergence and efficiency: the resulting designs should be certifiably optimal in comparison with all alternate designs at reduced computational cost for the optimization. The risk measures satisfying either of the certifiability conditions are classified under CRiBDO while satisfying both conditions makes the resulting optimal designs strongly certifiable. 
We compare and contrast the existing RBDO formulation based on probability of failure (PoF) with risk-based optimization formulations using the buffered probability of failure (bPoF) and the superquantile (a.k.a.\ conditional value-at-risk) risk measures. We show that RBDO does not satisfy either of the certifiability conditions while superquantiles and bPoF lead to CRiBDO formulations. An additional advantage of bPoF is the intuitive relation to PoF, which allows widely used PoF-based optimization formulations to be easily transitioned to bPoF-based optimization formulations.

Both bPoF and superquantile risk measures introduce data-informed conservativeness by encoding extra information about the limit state function distribution in the form of the magnitude of failure. This highlights a way to get desirable conservativeness in engineering risk-based optimization by switching to a different optimization formulation and getting certifiably risk-averse designs. Superquantiles and bPoF provide alternate measures of risk that avoid hard-threshold characterizations of failure events used in PoF, thus, relaxing the guesswork associated with picking limit state function thresholds. We highlight the data-informed nature of the conservativeness through numerical experiments on bPoF-based CRiBDO for short column structural design problem and superquantile-based CRiBDO for cooling fin thermal design problem.

When the underlying limit state and objective functions are convex w.r.t.\ the design variables, using bPoF and superquantiles can lead to strongly certifiable designs because these risk measures can preserve convexity and lead to certifiably optimal designs. We show this property for bPoF-based CRiBDO through convex reformulation of the limit state and the objective function for the short column design problem. The bPoF-based CRiBDO results in a convex optimization problem that can be solved using convex optimizers and provide global convergence guarantees. Although, such convex reformulations are not always possible, using convex approximations for the limit state functions is one way of addressing non-convex problems. 

Estimation of bPoF and superquantiles can be expensive, and recent results for Monte Carlo variance reduction techniques along with approximation techniques have been shown to be effective ways of reducing the computational effort for estimating superquantiles and bPoF \cite{JORoyset_LBonfiglio_GVernengo_SBrizzolara_2017a,bonfiglio2019multidisciplinary,MMHarajli_RTRockafellar_JORoyset_2015a,DPKouri_TMSurowiec_2016a,ZZou_DPKouri_WAquino_2018a,HKT2020_Adaptive_ROM_CVAR_estimation,HKTQ18CVaRROMS,HYang_MGunzburger_2016a,chaudhuri2020multifidelity}. More such research efforts will lead to flexibility in switching to appropriate and advantageous CRiBDO formulations for designing safe engineering systems.

\section*{Acknowledgement}
This work has been supported in part by the Air Force Office of Scientific Research (AFOSR) MURI on managing multiple information sources of multi-physics systems award numbers FA9550-15-1-0038 and FA9550-18-1-0023, and Air Force Center of Excellence on Multi-Fidelity Modeling of Rocket Combustor Dynamics award FA9550-17-1-0195. The fourth author acknowledges the support from Office of Naval Research under MIPR N0001420WX00519.


\bibliography{RiskOpt_bib}
\bibliographystyle{aiaa}

\end{document}